\documentclass[a4paper,12pt]{article}     

\usepackage{graphicx}		  
\usepackage{amsmath,amssymb}	

\newtheorem{lemma}{Lemma}[section]
\newtheorem{theorem}[lemma]{Theorem}
\newtheorem{proposition}[lemma]{Proposition}
\newtheorem{corollary}[lemma]{Corollary}
\newtheorem{example}[lemma]{Example}

\newtheorem{remark}[lemma]{Remark}

\def\authorfont{\footnotesize}

\def\ccode#1{\par
\vspace*{8pt}
{\authorfont{\leftskip18pt\rightskip\leftskip
\noindent #1\par}}\par}

\newenvironment{Proof}{
\hspace*{-9mm}
{ \it Proof.}}
{\hfill {$\square$}\vspace{1.5em}}

\begin{document}

\begin{center}{
{\Large 
 Distinguishing surface-links with 4-charts with 2 crossings and 8 black vertices}
\vspace{10pt}
\\ 
Teruo NAGASE and Akiko SHIMA\footnote{The second author is supported by JSPS KAKENHI Grant Number 21K03255.}
}
\end{center}

\begin{abstract}
Charts are oriented labeled graphs in a disk.
Any simple surface braid (2-dimensional braid) can be described by using a chart.
Also, a chart represents an oriented closed surface 
(called a surface-link)
embedded in 4-space.
In this paper, we investigate surface-links
by using charts.
In \cite{StI}, \cite{StII},
we gave an enumeration of the charts with two crossings. 
In particular, there are two classes for 4-charts
with 2 crossings and 8 black vertices.
The first class represents surface-links each of which is connected. The second class represents surface-links each of which is exactly two connected components. 
In this paper, by using quandle colorings,
we shall show that the charts in the second class 
represent different surface-links.
\end{abstract}

%
%
%
%

\ccode{2020 Mathematics Subject Classification. Primary 57K45,05C10; Secondary 57M15.}
\ccode{ {\it Key Words and Phrases}. surface-link, chart, C-move, crossing. }


\setcounter{section}{0}
\section{Introduction}

Charts are oriented labeled graphs in a disk (see  \cite{KnottedSurfaces},\cite{BraidBook}, Section~\ref{s:Prel}  for the precise definition of charts).
Let $D_1, D_2$ be 2-dimensional disks.
Any simple surface braid (2-dimensoinal braid) can be described 
by using a chart,
here a simple surface braid is a properly embedded surface
$S$ in the 4-dimensional disk $D_1\times D_2$ such that
a natural map $\pi:S\subset D_1\times D_2\to D_2$ is 
a simple branched covering map of $D_2$ and
the boundary $\partial S$ is a trivial closed braid in
the solid torus $D_1\times \partial D_2$
(see \cite{BraidThree}, \cite{CharacterKamada}, \cite[Chapter 14 and Chapter 18]{BraidBook}, Section~\ref{s:SurfaceBraids}).
Also, from a chart, 
we can construct a simple closed surface braid in 4-space ${\Bbb R}^4$. This surface is an oriented closed surface 
embedded in ${\Bbb R}^4$.
On the other hand, any oriented embedded closed surface 
 in ${\Bbb R}^4$ is ambient isotopic to a simple
closed surface braid
 (see \cite{CharacterKamada},\cite[Chapter 23]{BraidBook}). 
A C-move 
is a local modification between two charts
in a disk (see Section~\ref{s:Prel} for C-moves).
A C-move between two charts induces 
an ambient isotopy between oriented closed surfaces 
corresponding to the two charts.
In this paper, we investigate oriented closed surfaces in 4-space
by using charts.

We will work in the PL category or smooth category. All submanifolds are assumed to be locally flat.

A {\it surface-link} is a closed surface embedded in 4-space ${\Bbb R}^4$. An orientable surface-link is called a {\it ribbon surface-link} if there exists an immersion of a 3-manifold $M$ into ${\Bbb R}^4$ sending the boundary of $M$ onto the surface-link such that each connected component of $M$ is a handlebody and its singularity consists of ribbon singularities, here a ribbon singularity is a disk in the image of $M$ whose pre-image consists of two disks; one of the two disks is a proper disk of $M$ and the other is a disk in the interior of $M$. 
In the words of charts, a ribbon surface-link is a surface-link corresponding to a {\it ribbon chart}, a chart C-move equivalent to a chart without white vertices \cite{BraidThree}.
Here, in a chart there are 
three kinds of vertices; black vertices (vertices of degree 1), crossings (veritces of degree 4) and white vertices (veritces of degree 6).

Kamada showed that any 3-chart is a ribbon chart \cite{BraidThree}. Kamada's result was extended by Nagase and Hirota: Any 4-chart with at most one crossing is a ribbon chart \cite{NagaseHirota}. We showed that any $n$-chart with at most one crossing is a ribbon chart \cite{MinimalChartOneCrossing}.
A {\it $2$-link} is a surface-link each of whose connected component
is a 2-sphere.
We also showed that if a chart with at most two crossings represents a 2-link, then
the chart is a ribbon chart \cite{NSTwoCrossing}, \cite{NSTwoCrossingII}.
We investigated the structure of c-minimal charts with
two crossings 
(see Section~\ref{s:Prel} for the precise definition of a c-minimal chart),
and gave an enumeration of the charts with two crossings in
\cite{StI}, \cite{StII}. 

\begin{figure}
\begin{center}
\includegraphics{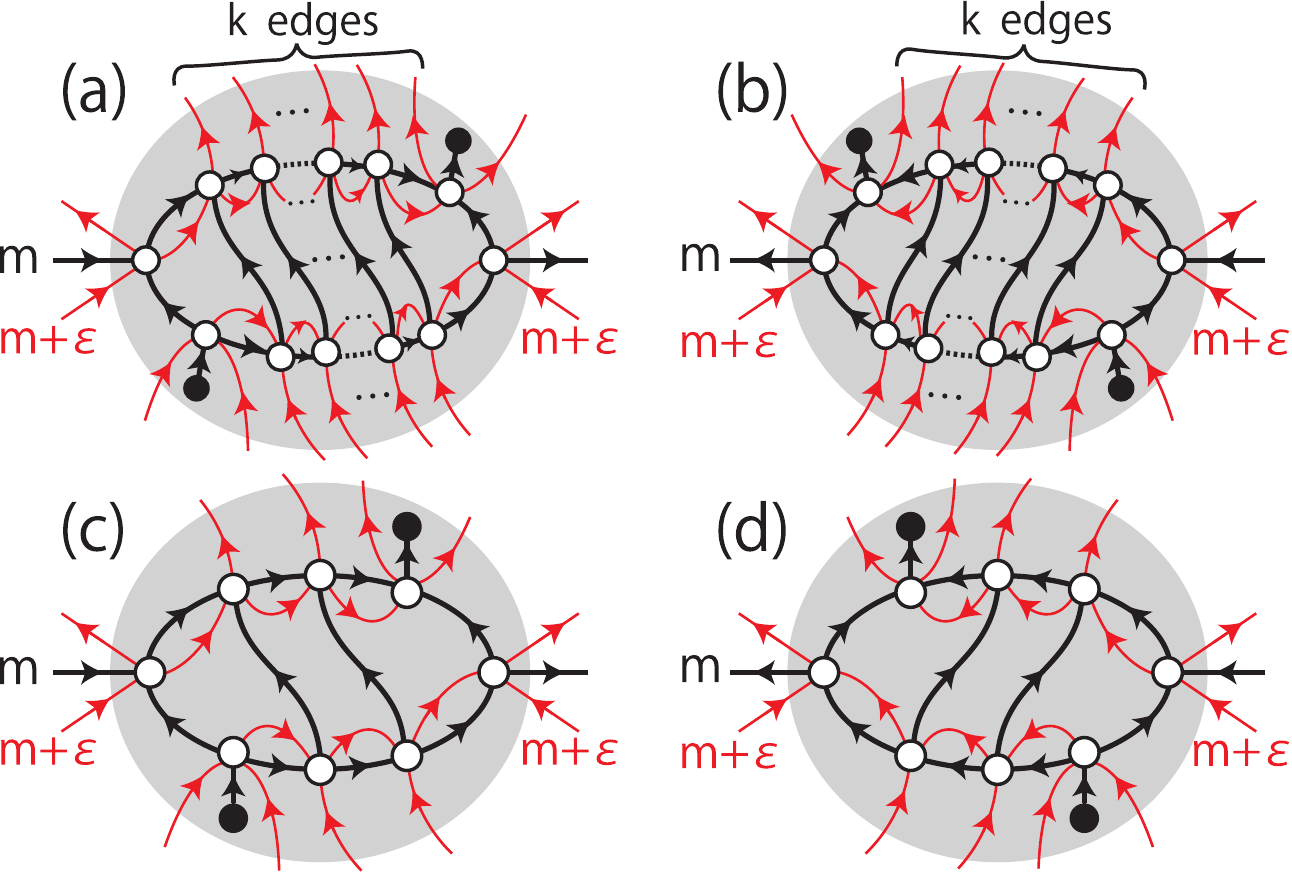}
\end{center}
\caption{ \label{Fig01} 
Type-I elementary IO-tangles of label $m$, here
$\varepsilon\in\{+1,-1\}$. 
(a),(b) Type-$I_k$.
(c),(d) Type-$I_3$.}
\end{figure}

Let $\Gamma$ be a chart in a 2-disk $D^2$, and $m$ a label of $\Gamma$.
Let $D$ be a 2-disk in $D^2$.
If the pair $(\Gamma\cap D,D)$ is one of the two pairs 
as shown in Fig.~\ref{Fig01}(a),(b), then the pair
is called a {\it Type-I elementary IO-tangle of label $m$}.
In particular,
if the disk $D$ contains exactly $2k+2$ white vertices,
then the pair is called a {\it Type-$I_k$ elementary IO-tangle}.
For example, the two pairs as shown in Fig.~\ref{Fig01} (c),(d) are Type-$I_3$ elementary IO-tangles of label $m$.

A {\it hoop} is a closed edge of a chart
that contains neither crossings nor white vertices.

Let $\Gamma$ be a chart.
If $\Gamma'$ is obtained from $\Gamma$ by introducing a hoop surrounding $\Gamma$,
then we say that $\Gamma'$ is obtained from $\Gamma$
by a {\it conjugation}.
Note that
a conjugation also 
induces an ambient isotopy between surface-links 
corresponding to the two charts.

\begin{figure}
\begin{center}
\includegraphics{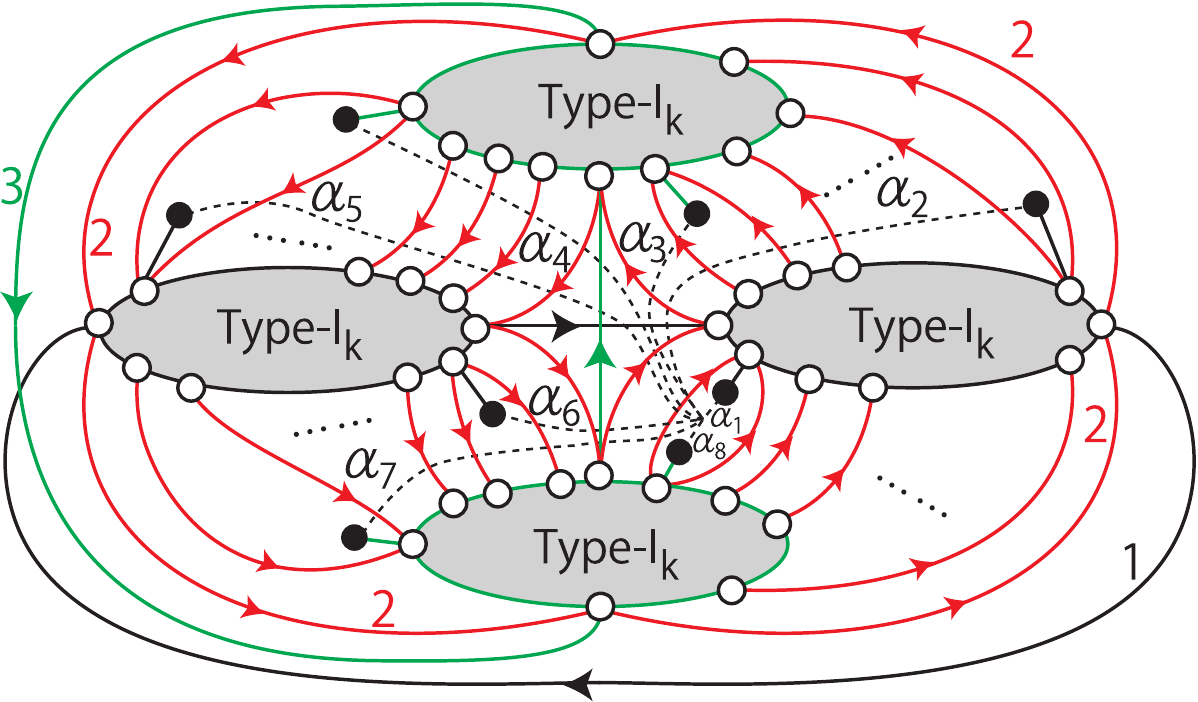}
\end{center}
\caption{ \label{Fig02} 
A 4-chart $T_k$ with 2 crossings and 8 black vertices. 
The dotted lines are a Hurwizt arc system 
$(\alpha_1,\alpha_2,\cdots,\alpha_8)$. }
\end{figure}

\begin{figure}
\begin{center}
\includegraphics{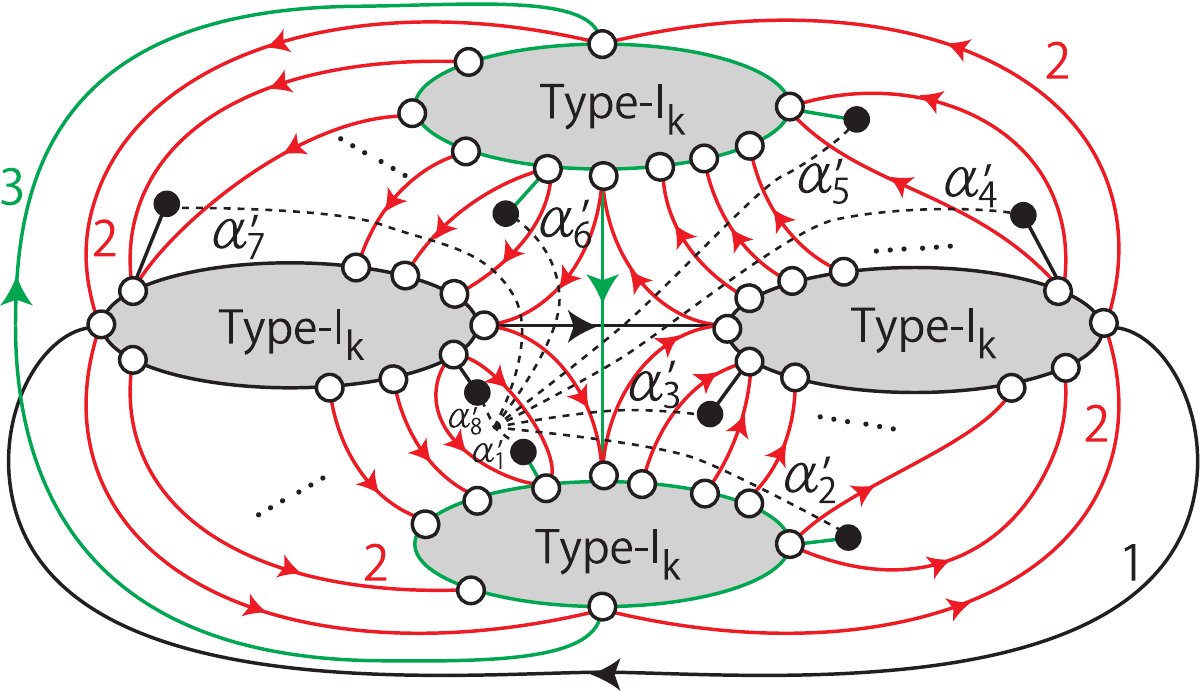}
\end{center}
\caption{ \label{Fig03} 
A 4-chart $T_k^*$ with 2 crossings and 8 black vertices. 
The dotted lines are a Hurwizt arc system 
$(\alpha_1',\alpha_2',\cdots,\alpha_8')$.}
\end{figure}

We denote by $T_k$ (resp. $T_k^*$)
a 4-chart as shown in Fig.~\ref{Fig02}
(resp. Fig.~\ref{Fig03}).
These 4-charts have 2 crossings and 8 black vertices,
and four Type-$I_k$ elementary IO-tangles.
The following is a special case in \cite[Section 9]{StII}.

\begin{figure}
\begin{center}
\includegraphics{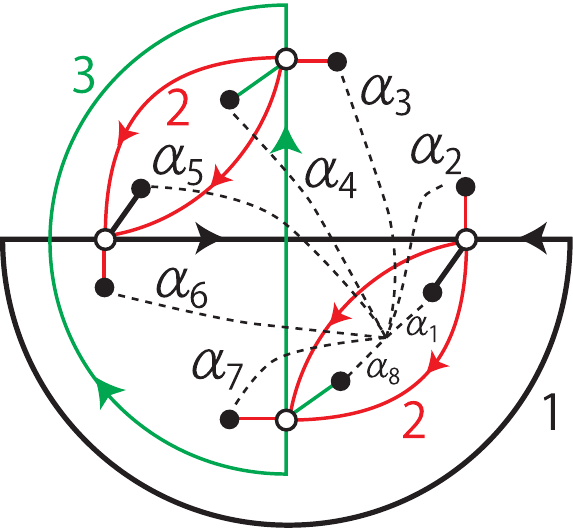}
\end{center}
\caption{ \label{Fig04} 
A 4-chart $T_0$ with 2 crossings and 8 black vertices.
The dotted lines are a Hurwizt arc system 
$(\alpha_1,\alpha_2,\cdots,\alpha_8)$. }
\end{figure}

\begin{theorem}
\label{CorTwoCrossing}
{\rm $($\cite[Section 9]{StII}$)$}
Let $\Gamma$ be a c-minimal $4$-chart with exactly
two crossings.
If $\Gamma$ contains exactly $8$ black vertices,
then $\Gamma$ is obtained from
$T_k$ or $T_k^*$ for some $k$,
or a $4$-chart as shown in Fig.~\ref{Fig04}
by C-moves and conjugations. 
\end{theorem}

We denote by $T_0$
the 4-chart as shown in Fig.~\ref{Fig04}.
The $4$-chart $T_0$ is a c-minimal chart (see {\rm \cite{NST}}).
However, we do not know that 4-charts $T_k,T_k^*$ are c-minimal charts. 

In this paper, by using quandle colorings of surface-links,
we shall show that
the 4-charts $T_0,T_1,T_3,T_5,\cdots,T_{2k-1},\cdots$
are not C-move equivalent each other,
and the 4-charts $T_0,T^*_1,T^*_3,T^*_5,\cdots,T^*_{2k-1},\cdots$
are not C-move equivalent each other 
(see Theorem~\ref{MainTheorem} and 
Corollary~\ref{CorMainTheorem}).

The paper is organized as follows.
In Section~\ref{s:Prel}, 
we define charts and c-minimal charts.
In Section~\ref{s:QuandleColor},
we review the definition of quandles and quandle colorings.
In Section~\ref{s:IsomorphismFreeQuandle}, 
for each classical $n$-braid $b$,
we shall define the quandle isomorphism
$Q(b)$ of the free quandle $F_Q\langle x_1,x_2,\cdots,x_n\rangle$
generated by $x_1,x_2,\cdots,x_n$. 
These maps $Q(b)$ are important to count 
the number of quandle colorings. 
In Section~\ref{s:SurfaceBraids},
we review surface braids and 
the construction of surface-links from charts.
In Section~\ref{s:KnotQuandle},
 we review knot quandles and a relation between
quandle colorings and knot quandles.
We shall calculate knot quandles of surface braids described 
by the charts $T_k$ and $T_k^*$.
In Section~\ref{s:T0},
we shall count the number of quandle colorings of a surface braid described by the chart $T_0$.
In Section~\ref{s:T2k},
we shall count the number of quandle colorings of a surface braid described by the chart $T_{2k}$.
In Section~\ref{s:T2k-1},
we shall count the number of quandle colorings of a surface braid described by the chart $T_{2k-1}$. We shall show the main theorem (Theorem~\ref{MainTheorem}).


\section{Preliminaries}
\label{s:Prel}

In this section, 
we introduce 
the definition of charts and its related words.

\begin{figure}
\begin{center}
\includegraphics{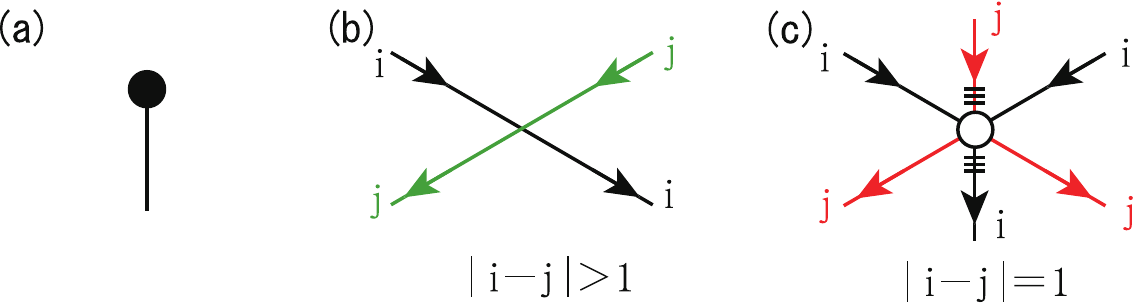}
\end{center}
\caption{ \label{Fig05} (a) A black vertex. (b) A crossing. (c) A white vertex. 
Each arc with three transversal short arcs is a middle arc at the white vertex. }
\end{figure}

Let $n$ be a positive integer.
An $n$-{\it chart}  
(a braid chart of degree $n$ \cite{KnottedSurfaces}
or a surface braid chart of degree $n$ \cite{BraidBook}) 
is 
an oriented labeled graph in the interior of a disk,
which may be empty 
or
have closed edges without vertices
satisfying the following four conditions
(see Fig.~\ref{Fig05}):
\begin{enumerate}
\item[(i)] 
Every vertex has degree $1$, $4$, or $6$.
\item[(ii)] 
The labels of edges are 
in $\{1,2,\dots,n-1\}$.
\item[(iii)]
In a small neighborhood of
each vertex of degree $6$,
there are six short arcs,
three consecutive arcs are
oriented inward 
and
the other three are outward,
and
these six are labeled $i$ and $i+1$
alternately for some $i$,
where the orientation and label of
each arc are inherited from
the edge containing the arc.
\item[(iv)]
For each vertex of degree $4$,
diagonal edges have the same label
and
are oriented coherently,
and the labels $i$ and $j$ of
the diagonals satisfy $|i-j|>1$.
\end{enumerate}
We call a vertex of degree $1$ a {\it black vertex},
a vertex of degree $4$ a {\it crossing}, and 
a vertex of degree $6$ a {\it white vertex}
respectively.

Among six short arcs
in a small neighborhood of
a white vertex,
a central arc of each three consecutive arcs
oriented inward (resp. outward) 
is called a   
{\it middle arc} at the white vertex
(see Fig.~\ref{Fig05}(c)).


\begin{figure}
\begin{center}
\includegraphics{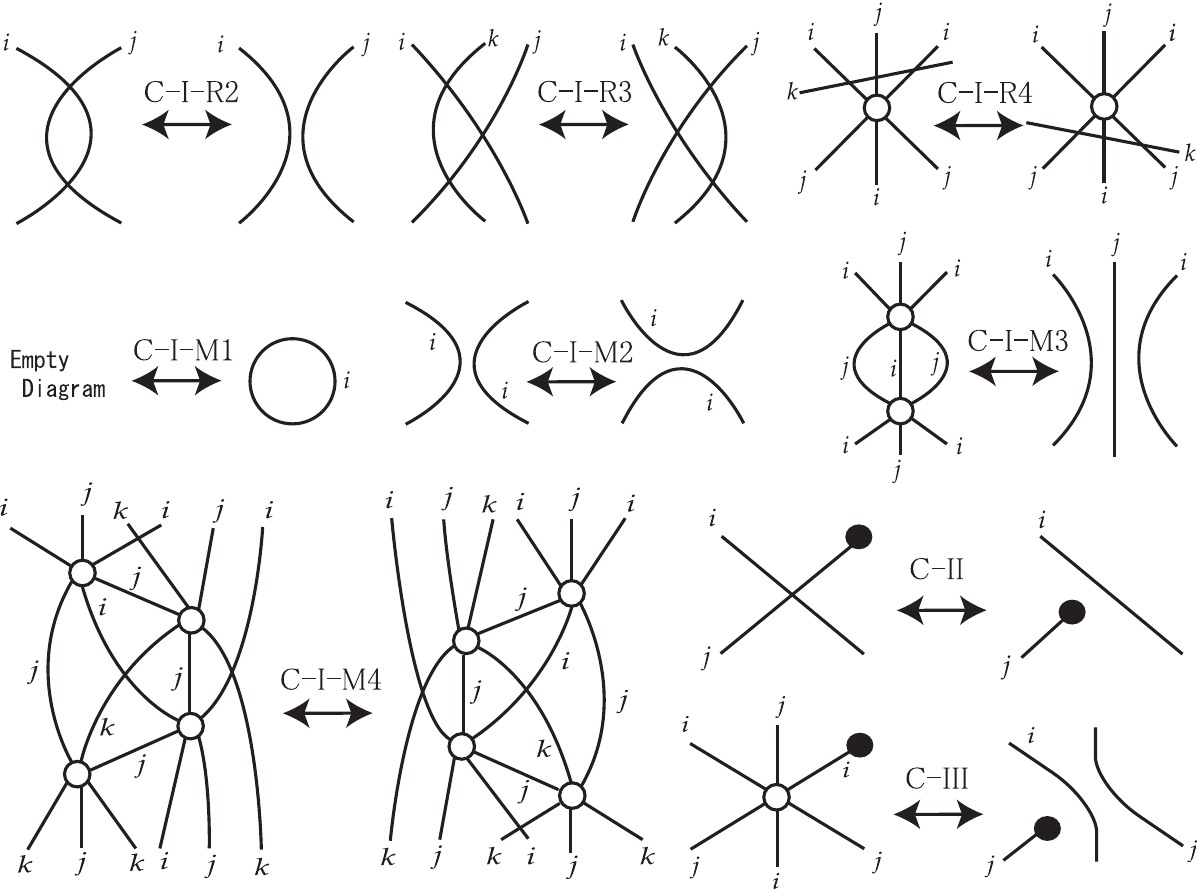}
\end{center}
\caption{ \label{Fig06} For the C-III move, 
the edge with the black vertex does not contain a middle arc at
a white vertex in the left figure. }
\end{figure}

Now {\it C-moves} are local modifications 
of charts as shown in Fig.~\ref{Fig06}
(cf. \cite{KnottedSurfaces}, 
\cite{BraidBook} and \cite{Tanaka}).
Two charts are said to be {\it C-move equivalent}  if there exists
a finite sequence of C-moves 
which modifies one of the two charts 
to the other.

\begin{figure}
\begin{center}
\includegraphics{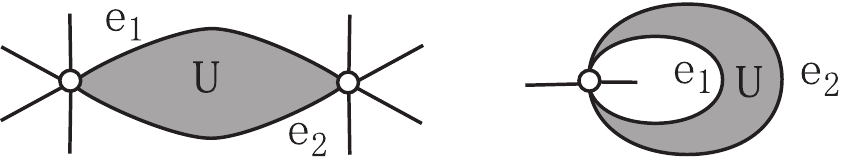}
\end{center}
\caption{ \label{Fig07} Bigons.}
\end{figure}

Let $\Gamma$ be a chart. 
Let $e_1$ and $e_2$ be edges of $\Gamma$
which connect two white vertices $w_1$ and $w_2$
where possibly $w_1=w_2$.
Suppose that 
the union $e_1\cup e_2$ bounds 
an open disk $U$.
Then $Cl(U)$ 
is called 
a {\it bigon} of $\Gamma$
provided that
any edge containing $w_1$ or $w_2$ 
does not intersect the open disk $U$
(see Fig.~\ref{Fig07}).
Note that neither $e_1$ nor $e_2$ contains a crossing.

An edge in a chart is called 
a {\it free edge}
if it has
two black vertices.

Let $\Gamma$ be a chart.
Let 
$c(\Gamma)$, $w(\Gamma)$, $f(\Gamma)$, 
and 
$b(\Gamma)$ be 
the number of crossings, 
the number of white vertices, 
the number of free edges, 
and 
the number of bigons of $\Gamma$
respectively.
The 4-tuple $(c(\Gamma),w(\Gamma),-f(\Gamma),-b(\Gamma))$ is called a 
{\it $c$-complexity} of the chart $\Gamma$.
A chart $\Gamma$ is called a {\it c-minimal chart} if its complexity is minimal among the charts C-move equivalent to the chart $\Gamma$ with respect to the lexicographic order of 4-tuples of integers.

The following is a special case in \cite[Section 9]{StII}.
From the following theorem,
we have Theorem~\ref{CorTwoCrossing}.

\begin{theorem}{\rm $($\cite[Section 9]{StII}$)$}
Let $\Gamma$ be a c-minimal $4$-chart with exactly
two crossings.
If $\Gamma$ contains exactly $8$ black vertices,
then $\Gamma$ is C-move equivalent to
a union of hoops and one of the following charts:
\begin{enumerate}
\item[{\rm (a)}] a $4$-chart $T_k$ for some $k$ $($see Fig.~\ref{Fig02}$)$,
\item[{\rm (b)}] a $4$-chart $T_k^*$ for some $k$ $($see Fig.~\ref{Fig03}$)$, or 
\item[{\rm (c)}] a $4$-chart $T_0$ as shown in Fig.~\ref{Fig04}. 
\end{enumerate}
\end{theorem}

In this paper 
we shall show that
the 4-charts $T_0,T_1,T_3,T_5,\cdots,T_{2k-1},\cdots$
are not C-move equivalent each other,
and the 4-charts $T_0,T^*_1,T^*_3,T^*_5,\cdots,T^*_{2k-1},\cdots$
are not C-move equivalent  each other.

In this paper
for a subset $X$ in a space
we denote 
the closure of $X$
by $Cl(X)$.
For a manifold $M$,
we denote the boundary of $M$
by $\partial M$.

\section{Quandle colorings}
\label{s:QuandleColor}

In this section, we review the definition of quandles and quandle colorings. We shall introduce the main theorem (Theorem~\ref{MainTheorem}) in this paper.

A {\it quandle} is defined to be a set $Q$ with 
a binary operation $*:Q\times Q\to Q$ such that
\begin{enumerate}
\item[(i)] for any $a\in Q$, we have $a*a=a$, 
\item[(ii)] for any $a,b\in Q$, 
there is a unique $x\in Q$ such that $x*a=b$, and
\item[(iii)] for any $a,b,c\in Q$, we have $(a*b)*c=(a*c)*(b*c)$ .
\end{enumerate}
By Condition (ii),
there exists unique element $x$ with $x*a=b$
for $a,b\in Q$.
We denote $x=b\overline{*}a$.

\begin{example}
\label{PseudoTrivialQuandle}
{\rm Let $N$ be an integer with $N\geqq3$.
Let $Q_N=\{1,2,\cdots,N\}$.
The binary operation $*$ is defined by
$$\left\{\begin{array}{lcl}
x*y =x, && ({\rm if} \ y\not=N),\\
x*N =x+1, && ({\rm if} \ x\not=N-1,x\not=N),\\
(N-1)*N =1, && \\
N*N =N. && 
\end{array}
\right.$$
Then $(Q_N,*)$ is a quandle. 
Note that $N*y=N$ for all $y\in Q_N$
and 
$$\left\{\begin{array}{lcl}
x\overline{*}y =x, && ({\rm if} \ y\not=N),\\
x\overline{*}N =x-1, && ({\rm if} \ x\not=1,x\not=N),\\
1\overline{*}N =N-1, && \\
N\overline{*}N =N. && 
\end{array}
\right.$$
 For example, if $N=4$,
then the binary operation $*$ is as follows (see Table~1).}
\end{example}

\begin{center}
\begin{tabular}{|c|cccc|}
\hline
 & 1 & 2 & 3 & 4\\
\hline
  1 & 1 & 1 & 1 & 2\\
  2 & 2 & 2 & 2 & 3\\
  3 & 3 & 3 & 3 & 1\\
  4 & 4 & 4 & 4 & 4\\
\hline
\end{tabular}\vspace{2mm}

Table 1: An example $(Q_4,*)$ of a quandle. 
\end{center}

\begin{figure}
\begin{center}
\includegraphics{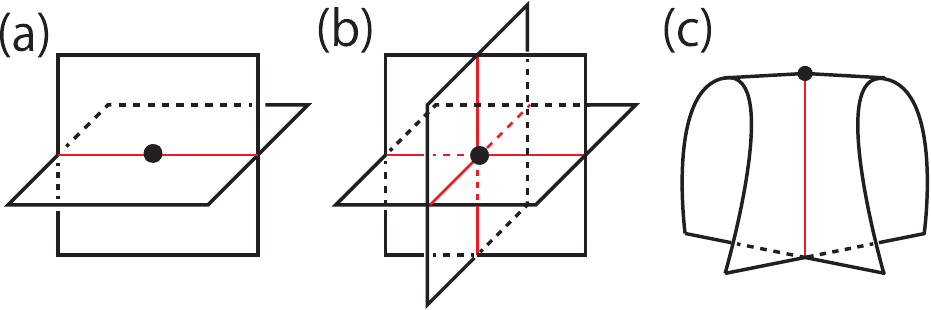}
\end{center}
\caption{ \label{Fig08} 
(a) a double point, (b) a triple point, (c) a branch point. }
\end{figure}

A {\it surface-link} is a connected or disconnected closed 
 surface embedded in 4-space ${\Bbb R}^4$
locally flatly.
In particular, a surface-link is {\it oriented} if
each connected component is an oriented closed surface.
A projection $\pi:{\Bbb R}^4\to{\Bbb R}^3$
is {\it generic} for a surface-link $F$ if
the image $\pi(F)$ is locally homeomorphic to
(i) a single sheet,
(ii) two transversely sheets,
(iii) three transversely sheets, or
(iv) a cross-cup.
The points corresponding to (ii), (iii) and (iv)
are called a {\it double point}, a {\it triple point},
and a {\it branch point} of the generic projection,
 respectively. 
See Fig.~\ref{Fig08}.
A {\it surface diagram} of $F$ is the image $\pi(F)$
with additional crossing information at the singularity set.
There are two intersecting sheets near a double point,
one of which is higher than the other with respect to $\pi$.
They are called an {\it over-sheet} and an {\it under-sheet},
respectively.
In order to indicate crossing information,
we break the under sheet into two pieces missing the over-sheet.
Then the surface diagram is presented by disjoint union of compact surfaces which are called {\it broken sheets}
(see \cite[page 2]{KnottedSurfaces}, \cite[page 58]{BraidBook}).

\begin{figure}
\begin{center}
\includegraphics{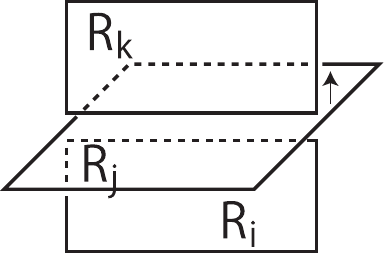}
\end{center}
\caption{ \label{Fig09} 
A quandle relation ${\mathcal C}(R_i)*{\mathcal C}(R_j)={\mathcal C}(R_k)$ at a double point. }
\end{figure}

Let $Q$ be a fixed quandle.
Let $D$ be a surface diagram of an oriented surface-link.
Let ${\mathcal R}$ be the set of all broken sheets of $D$.
The normals are given in such a way that
(tangents, normal) matches the orientation of 3-space ${\Bbb R}^3$, see Fig.~\ref{Fig09}.
A {\it quandle coloring} ${\mathcal C}$
is a map ${\mathcal C}:{\mathcal R}\to Q$ such that 
at each double point,
there exist three broken sheets $R_i,R_j,R_k\in {\mathcal R}$ uniquely such that 
$${\mathcal C}(R_i)*{\mathcal C}(R_j)={\mathcal C}(R_k)$$
where
$R_j$ is the over-sheet at the double point,
$R_i,R_k$ are under-sheets at the double point,
the normal of $R_j$ points from $R_i$ to $R_k$
(see Fig.~\ref{Fig09}).

Let $Q$ be a quandle, $F$ an oriented surface-link.
Let  ${\rm Col}_Q(F)$ be the set of colorings of the surface-link diagram of $F$.
Then the cardinality $|{\rm Col}_Q(F)|$ of colorings is
a surface-link invariant (see \cite[Page 131,132]{Surfaces4Space}).
Note that if an oriented surface-link $F$ is a trivial surface-link with $k$-components,
then we have $|{\rm Col}_Q(F)|=N^k$ for any quandle $Q$ with
exactly $N$ elements.
Here a {\it trivial} surface-link means
the surface-link bounds a disjoint union of handlebodies 
in 4-space. 

Let $\Gamma$ be a chart.
Then we can construct a surface-link embedded in 4-space 
from $\Gamma$. In this paper we denote it by $F(\Gamma)$ 
(a closure of a surface braid described by a chart $\Gamma$,
see Section~\ref{s:SurfaceBraids}). 

The following is the main result in this paper.

\begin{theorem}
\label{MainTheorem} 
Let $N,k$ be integers with $N\geqq 3$ and $k\geqq1$.
Let $T_0,T_k,T_k^*$ be $4$-charts as shown in 
Fig.~\ref{Fig02}, Fig.~\ref{Fig03}, Fig.~\ref{Fig04},
respectively.
For the quandles $Q_N,Q_{k+2}$ in Example~\ref{PseudoTrivialQuandle},
we have the following:

\begin{enumerate}
\item[{\rm (a)}] $|{\rm Col}_{Q_N}(F(T_0))|=(N-1)^2+1$. 
\item[{\rm (b)}] $|{\rm Col}_{Q_N}(F(T_{2k}))|=N$,
$|{\rm Col}_{Q_N}(F(T_{2k}^*))|=N$.
\item[{\rm (c)}] 
$|{\rm Col}_{Q_{k+2}}(F(T_{2k-1}))|=(k+2)^2$,
$|{\rm Col}_{Q_{k+2}}(F(T_{2k-1}^*))|=(k+2)^2$.
\item[{\rm (d)}] 
$|{\rm Col}_{Q_{k+2}}(F(T_{2\ell-1}))|=(k+1)^2+1$,
$|{\rm Col}_{Q_{k+2}}(F(T_{2\ell-1}^*))|=(k+1)^2+1$ for any integer $\ell$ with $1\leqq \ell <k$. 
\end{enumerate}
\end{theorem}

Hence we have the following corollary from the above theorem.

\begin{corollary}
\label{CorMainTheorem}
For the $4$-charts $T_0,T_{2k-1},T_{2k-1}^*$,
we have the following:
 
\begin{enumerate}
\item[{\rm (a)}]
The surface-links $F(T_0),F(T_1),F(T_3),F(T_5),\cdots,F(T_{2k-1}),\cdots$ are not ambient isotopic each other.
\item[{\rm (b)}]
The surface-links $F(T_0),F(T_1^*),F(T_3^*),F(T_5^*),\cdots,F(T_{2k-1}^*),\cdots$ are not ambient isotopic each other.
\item[{\rm (c)}]
The $4$-charts $T_0,T_1,T_3,T_5,\cdots,T_{2k-1},\cdots$ are not C-move equivalent each other.
\item[{\rm (d)}]
The $4$-charts $T_0,T_1^*,T_3^*,T_5^*,\cdots,T_{2k-1}^*,\cdots$ are not C-move equivalent each other.
\end{enumerate}
\end{corollary} 

\begin{remark}
\begin{enumerate}
\item[{\rm (1)}]
Each of the surface-links $F(T_0),F(T_{2k-1}),F(T_{2k-1}^*)$ 
is exactly two connected components.

\item[{\rm (2)}] Each of the surface-links $F(T_{2k}),F(T_{2k}^*)$ is connected.

\item[{\rm (3)}] The $4$-chart $T_0$ is a c-minimal chart
$($see {\rm \cite{NST}}$)$.

\item[{\rm (4)}] We do not know that $T_{2k-1},T_{2k-1}^*,T_{2k},T_{2k}^*$ are c-minimal charts. 
\end{enumerate}
\end{remark}

\section{Isomorphisms of the free quandle}
\label{s:IsomorphismFreeQuandle}

In this section, for each classical $n$-braid $b$,
we shall define the quandle isomorphism
$Q(b)$ of the free quandle $F_Q\langle x_1,x_2,\cdots,x_n\rangle$ (see \cite[Example 10.8]{QuandleInvariant}).

Let $Q,Q'$ be quandles.
If a map $f:Q\to Q'$ satisfies $f(x*y)=f(x)*f(y)$ for all $x,y\in Q$,
then $f$ is called a {\it homomorphism}.
Moreover, if $f$ is bijective,
then $f$ is called an {\it isomorphism}.

Let $B_n$ be the $n$-braid group generated 
by standard generators $\sigma_1,\sigma_2,\cdots,\sigma_{n-1}$
whose relations are 
$\sigma_i\sigma_j\sigma_i=\sigma_j\sigma_i\sigma_j$ for
$|i-j|=1$ and
$\sigma_i\sigma_j=\sigma_j\sigma_i$ for
$|i-j|>1$.

Let $F_Q\langle x_1,x_2,\cdots,x_n\rangle$ be the free quandle
generated by $x_1,x_2,\cdots,x_n$.
For the trivial braid $1$ of the $n$-barid group $B_n$,
we define the map $Q(1):F_Q\langle x_1,x_2,\cdots,x_n\rangle\to F_Q\langle x_1,x_2,\cdots,x_n\rangle$
by $Q(1)(x_i)=x_i$ for all $i$
(i.e. $Q(1)$ is the identity map).

By inductive,
we define the isomorphism of $F_Q\langle x_1,x_2,\cdots,x_n\rangle$ for a braid as follows:
For a braid $b$ in $B_n$ and 
 a generator of $\sigma_i$ of $B_n$,
if we have an isomorphism $Q(b)$ of $F_Q\langle x_1,x_2,\cdots,x_n\rangle$,
then we define the isomorphisms $Q(\sigma_i b)$ and $Q(\sigma_i^{-1}b)$ as follows.
Let $(a_1,a_2,\cdots,a_n)=(Q(b)(x_1),Q(b)(x_2),\cdots,Q(b)(x_n))$.
Then 
$$Q(\sigma_i b)(x_j)=
\left\{\begin{array}{lcl}
 a_{i+1}\overline{*}a_i, & & ({\rm if} \ j=i),\\
a_i, & & ({\rm if} \ j=i+1),\\ 
a_j,  & & ({\rm otherwise}).\\
\end{array}\right.$$
$$Q(\sigma_i^{-1}b)(x_j)=
\left\{\begin{array}{lcl}
a_{i+1}, & & ({\rm if} \ j=i),\\
a_i* a_{i+1}, & & ({\rm if} \ j=i+1),\\ 
a_j,  & & ({\rm otherwise}).\\
\end{array}\right.$$
Note that if
$b$ and $b'$ are the same $n$-braids,
then $Q(b)$ and $Q(b')$ are the same isomorphisms.

\begin{figure}
\begin{center}
\includegraphics{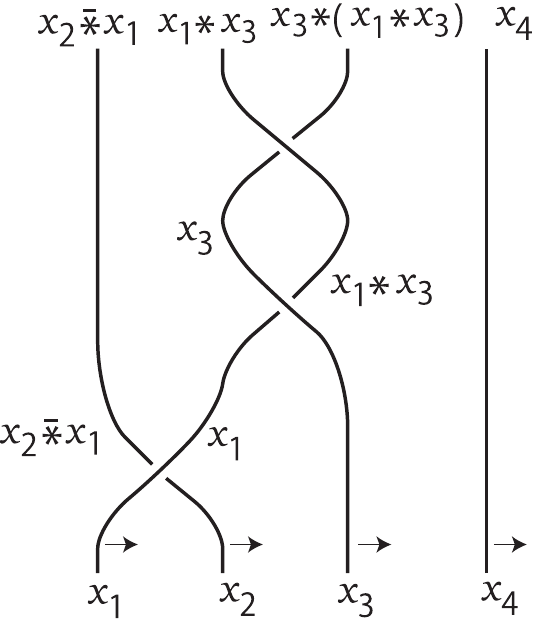}
\end{center}
\caption{ \label{Fig10} 
The 4-braid  $b=\sigma_2^{-2}\sigma_1$ in $B_4$.}
\end{figure}

For example, if $b=\sigma_2^{-2}\sigma_1\in B_4$,
then the quandle isomorphism $Q(b)$ maps
the generators $x_1,x_2,x_3,x_4$ as follows:
$$Q(b)(x_1)=x_2\overline{*}x_1, \ 
Q(b)(x_2)=x_1*x_3, \ 
Q(b)(x_3)=x_3 *(x_1*x_3), \ Q(b)(x_4)=x_4.$$
This map $Q(b)$ is obtained from the coloring of the braid
$b$ as shown in Fig.~\ref{Fig10}.

\begin{lemma}
\label{QNiPi+1N}
Let $N$ be an integer with $N\geqq 3$.
Let $Q_N$ be the quandle in Example~\ref{PseudoTrivialQuandle}, and
$c:F_Q\langle x_1,x_2,\cdots,x_n\rangle \to Q_N$ 
a homomorphism.
Let $b$ be an $n$-braid in $B_n$.
If $c(Q(b)(x_i))=1$ and $c(Q(b)(x_{i+1}))=N$,
then we have the following:

\begin{enumerate}
\item[{\rm (a)}]
$c(Q(\sigma_i^{-2k}b)(x_i))=\left\{
\begin{array}{ll}
k+1,  & ({\rm if} \ k=0,1,2,\cdots,N-2),\\
1, \ & ({\rm if} \ k=N-1).\\
\end{array}\right.$
\item[{\rm (b)}]
$c(Q(\sigma_i^{-2k}b)(x_{i+1}))=N$ for $k=0,1,\cdots,N-1$.
\item[{\rm (c)}]
$c(Q(\sigma_i^{-(2k+1)}b)(x_{i}))=N$ for $k=0,1,\cdots,N-2$.
\item[{\rm (d)}]
$c(Q(\sigma_i^{-(2k+1)}b)(x_{i+1}))=\left\{
\begin{array}{ll}
k+2,  & ({\rm if} \ k=0,1,2,\cdots,N-3),\\
1, \ & ({\rm if} \ k=N-2).\\
\end{array}\right.$
\end{enumerate}
\end{lemma}

\begin{Proof}
We can easily check for the case $k=0$.

Now suppose (a) and (b) hold  for the case $k=s$.
Then by (b)
$$c(Q(\sigma_i^{-(2s+1)}b)(x_i)) = c(Q(\sigma_i^{-2s}b)(x_{i+1}))=N.$$
Thus (c) holds for the case $k=s$.

If $s\leqq N-3$, then $s+1\leqq N-2$. 
Thus by (a) and (b) we have\vspace{2mm}

$\begin{array}{rcl}
c(Q(\sigma_i^{-(2s+1)}b)(x_{i+1}))& =& c(Q(\sigma_i^{-2s}b)(x_i)*Q(\sigma_i^{-2s}b)(x_{i+1}))\vspace{2mm}\\
& = & (s+1)*N=s+2.
\end{array}$\vspace{2mm}\\
If $s=N-2$,
then by (a) and (b) we have\vspace{2mm}

$\begin{array}{rcl}
c(Q(\sigma_i^{-(2s+1)}b)(x_{i+1}))& =& c(Q(\sigma_i^{-2s}b)(x_i)*Q(\sigma_i^{-2s}b)(x_{i+1}))\vspace{2mm}\\
& = & (s+1)*N=(N-1)*N\vspace{2mm}\\
& = & 1.
\end{array}$\vspace{2mm}\\
Thus (d) holds for the case $k=s$.

Now suppose (c) and (d) hold  for the case $k=s$.
By (c), we have $c(Q(\sigma_i^{-(2s+1)}b)(x_i))=N$.
Thus \vspace{2mm}

$\begin{array}{rcl}
c(Q(\sigma_i^{-(2s+2)}b)(x_{i+1})) & = & c(Q(\sigma_i^{-(2s+1)}b)(x_i)*Q(\sigma_i^{-(2s+1)}b)(x_{i+1}))\\
 &= &  N*c(Q(\sigma_i^{-(2s+1)}b)(x_{i+1}))\\
& = & N.
\end{array}$\vspace{2mm}\\
Thus (b) holds for the case $k=s+1$.

By (d), we have

$\begin{array}{rcl}
c(Q(\sigma_i^{-(2s+2)}b)(x_i)) & = &
 c(Q(\sigma_i^{-(2s+1)}b)(x_{i+1}))\vspace{2mm}\\
& = & \left\{
\begin{array}{ll}
s+2,  & ({\rm if} \ s \leqq N-3),\\
1, \ & ({\rm if} \ s=N-2).\\
\end{array}\right.
\end{array}$\\
Hence (a) holds for the case $k=s+1$.
\end{Proof}

From the above lemma,
we have the following corollary:

\begin{corollary}
\label{CorQNiPi+1N}
Let $N,k$ be integers with $N\geqq 3$ and $k\geqq0$.
Let $Q_N$ be the quandle in Example~\ref{PseudoTrivialQuandle}, and
$c:F_Q\langle x_1,x_2,\cdots,x_n\rangle \to Q_N$ 
a homomorphism.
Let $b$ be an $n$-braid in $B_n$.
If $c(Q(b)(x_i))=1$ and $c(Q(b)(x_{i+1}))=N$,
then we have the following:

\begin{enumerate}
\item[{\rm (a)}]
If $k=s+m(N-1)$ $(0\leqq s\leqq N-2, s,m\in{\Bbb Z})$,
then $c(Q(\sigma_i^{-2k}b)(x_i))=s+1$.
\item[{\rm (b)}]
$c(Q(\sigma_i^{-2k}b)(x_{i+1}))=N$ for any $k$.
\item[{\rm (c)}]
$c(Q(\sigma_i^{-(2k+1)}b)(x_{i}))=N$ for any $k$.
\item[{\rm (d)}]
If $k=s+m(N-1)$ $(0\leqq s\leqq N-2, s,m\in{\Bbb Z})$, then

$c(Q(\sigma_i^{-(2k+1)}b)(x_{i+1}))=\left\{
\begin{array}{ll}
s+2,  & ({\rm if} \ s\not=N-2),\\
1, \ & ({\rm if} \ s=N-2).\\
\end{array}\right.$
\end{enumerate}
\end{corollary}

\begin{corollary}
\label{CorQNiNi+1P}
Let $N,k$ be integers with $N\geqq 3$ and $k\geqq0$.
Let $Q_N$ be the quandle in Example~\ref{PseudoTrivialQuandle}, and
$c:F_Q\langle x_1,x_2,\cdots,x_n\rangle \to Q_N$ 
a homomorphism.
Let $b$ be an $n$-braid in $B_n$.
If $c(Q(b)(x_i))=N$ and $c(Q(b)(x_{i+1}))=1$,
then we have the following:

\begin{enumerate}
\item[{\rm (a)}]
$c(Q(\sigma_i^{-2k}b)(x_i))=N$ for any $k$.
\item[{\rm (b)}]
If $k=s+m(N-1)$ $(0\leqq s\leqq N-2, s,m\in{\Bbb Z})$,
then $c(Q(\sigma_i^{-2k}b)(x_{i+1}))=s+1$.
\item[{\rm (c)}]
If $k=s+m(N-1)$ $(0\leqq s\leqq N-2, s,m\in{\Bbb Z})$,
then $c(Q(\sigma_i^{-(2k+1)}b)(x_i))=s+1$.
\item[{\rm (d)}]
$c(Q(\sigma_i^{-(2k+1)}b)(x_{i+1}))=N$ for any $k$.
\end{enumerate}
\end{corollary}

\begin{Proof}
We can easily check for $k=0$.
In particular, 
$c(Q(\sigma_i^{-1}b)(x_i))=c(Q(b)(x_{i+1}))=1$
and 
$c(Q(\sigma_i^{-1}b)(x_{i+1}))=c(Q(b)(x_i)*Q(b)(x_{i+1}))=N*1=N$.
Thus we can apply Corollary~\ref{CorQNiPi+1N}.

For any positive integer $k$,
by Corollary~\ref{CorQNiPi+1N}(b),(c), we have\vspace{2mm}

$\begin{array}{l}
c(Q(\sigma_i^{-2k}b)(x_i))=c(Q(\sigma_i^{-(2(k-1)+1)}\sigma_i^{-1}b)(x_i))=N,\\
c(Q(\sigma_i^{-(2k+1)}b)(x_{i+1}))=c(Q(\sigma_i^{-2k}\sigma_i^{-1}b)(x_{i+1}))=N.
\end{array}$\vspace{2mm}\\
Hence Statements~(a),(d) hold.

For a positive integer $k$, 
there exist integers $s,m$ with
$k=s+m(N-1)$ and $0\leqq s\leqq N-2$.
By Corollary~\ref{CorQNiPi+1N}(a), we have
$$c(Q(\sigma_i^{-(2k+1)}b)(x_i))=c(Q(\sigma_i^{-2k}\sigma_i^{-1}b)(x_i))=s+1.$$
Hence Statement~(c) holds.

If $s=0$, then $m\geqq1$ and $k-1=(N-2)+(m-1)(N-1).$
Thus by Corollary~\ref{CorQNiPi+1N}(d),
$$c(Q(\sigma_i^{-2k}b)(x_{i+1}))=c(Q(\sigma_i^{-(2(k-1)+1)}\sigma_i^{-1}b)(x_{i+1}))=1=0+1=s+1.$$
If $1\leqq s\leqq N-2$, then $k-1=(s-1)+m(N-1).$
Thus by Corollary~\ref{CorQNiPi+1N}(d),
$$c(Q(\sigma_i^{-2k}b)(x_{i+1}))=c(Q(\sigma_i^{-(2(k-1)+1)}\sigma_i^{-1}b)(x_{i+1}))=(s-1)+2=s+1.$$
Hence Statement~(b) holds.
\end{Proof}

\begin{lemma}
\label{QNii+1NotN}
Let $N$ be an integer with $N\geqq 3$.
Let $Q_N$ be the quandle in Example~\ref{PseudoTrivialQuandle}, and
$c:F_Q\langle x_1,x_2,\cdots,x_n\rangle \to Q_N$ 
a homomorphism.
Let $b$ be an $n$-braid in $B_n$.
If $c(Q(b)(x_i))\not=N$ and $c(Q(b)(x_{i+1}))\not=N$,
then
$$
(c(Q(\sigma_i^{\varepsilon} b)(x_i)),
c(Q(\sigma_i^{\varepsilon} b)(x_{i+1})))
=(c(Q(b)(x_{i+1})),c(Q(b)(x_i))),
$$
where $\varepsilon\in\{+1,-1\}$.
\end{lemma}

\begin{Proof}
We only show the case $\varepsilon=+1$.
By the definition of $Q(\sigma_i b)$, we have\vspace{2mm}\\
$\begin{array}{rcl}
(c(Q(\sigma_i b)(x_i)),c(Q(\sigma_i b)(x_{i+1}))) 
& = & 
(c(Q(b)(x_{i+1}))\overline{*}c(Q(b)(x_i)),c(Q(b)(x_i)))\vspace{2mm}\\
& = & (c(Q(b)(x_{i+1})),c(Q(b)(x_i))).
\end{array}$
\end{Proof}

By Lemma~\ref{QNii+1NotN}, we can show the following lemma:

\begin{lemma}
\label{QNi1i+12}
Let $N$ be an integer with $N\geqq 3$.
Let $Q_N$ be the quandle in Example~\ref{PseudoTrivialQuandle},
and $s,t$ elements in $Q_N$ different from $N$.
Let $c:F_Q\langle x_1,x_2,\cdots,x_n\rangle \to Q_N$ 
a homomorphism.
Let $b$ be an $n$-braid in $B_n$.
If $c(Q(b)(x_i))=s$ and $c(Q(b)(x_{i+1}))=t$,
then we have the following:

\begin{enumerate}
\item[{\rm (a)}]
$c(Q(\sigma_i^{2k}b)(x_i))=s$ and $c(Q(\sigma_i^{2k}b)(x_{i+1}))=t$ for any integer $k$.
\item[{\rm (b)}]
$c(Q(\sigma_i^{2k+1}b)(x_{i}))=t$ and $c(Q(\sigma_i^{2k+1}b)(x_{i+1}))=s$
for for any integer $k$.
\end{enumerate}
\end{lemma}

By the similar way of Lemma~\ref{QNii+1NotN},
we can show the following lemma:

\begin{lemma}
\label{QNii+1Equal}
Let $Q$ be a quandle, and
$c:F_Q\langle x_1,x_2,\cdots,x_n\rangle \to Q$ 
a homomorphism.
Let $b$ be an $n$-braid in $B_n$.
If $c(Q(b)(x_i))=c(Q(b)(x_{i+1}))$,
then
$$
(c(Q(\sigma_i^{\varepsilon} b)(x_i)),
c(Q(\sigma_i^{\varepsilon} b)(x_{i+1})))
=(c(Q(b)(x_{i})),c(Q(b)(x_{i+1}))),
$$
where $\varepsilon\in\{+1,-1\}$.
\end{lemma}

\section{Surface braids}
\label{s:SurfaceBraids}

In this section,
we review surface braids and 
the construction of surface-links from charts.

Let $D_1,D_2$ be 2-dimensional disks,
and $pr_i: D_1\times D_2\to D_i$
$(i=1,2)$ the $i$th factor projection.
Let $X_n$ be a set of $n$ interior points in $D_1$.
A {\it simple surface braid of degree $n$}
is an oriented 2-manifold $S$ embedded properly and locally flatly in $D_1\times D_2$
such that the restriction map $pr_2 |_S :S\to D_2$
is a simple branched covering map of degree $n$
and $\partial S=X_n\times \partial D_2$.
Here a branched covering map $f:S\to D_2$
of degree $n$ is {\it simple}
if $f^{-1}(y)$ consists of $n-1$ or $n$ points for each
point $y$ in $D_2$.

A surface braid $S$ of degree $n$
is extended to a closed surface $\widehat{S}$
in $D_1\times S^2$
such that $\widehat{S}\cap(D_1\times D_2)=S$
and $\widehat{S}\cap (D_1\times \overline{D}_2)=X_n\times \overline{D}_2$,
where $S^2$ is the 2-sphere obtained from $D_2$ by
attaching a 2-disk $\overline{D}_2$ along the boundary.
By identifying $D_1\times S^2$ with a regular neighborhood of a standard 2-sphere in 4-space ${\Bbb R}^4$,
we assume that $\widehat{S}$ is a closed oriented surface
embedded in ${\Bbb R}^4$.
We call it the {\it closure} of $S$ in ${\Bbb R}^4$.
It is proved in \cite{CharacterKamada}
that every closed oriented surface embedded in ${\Bbb R}^4$
is ambient isotopic to the closure of a surface braid.

Let $\Gamma$ be an $n$-chart in a 2-disk $D_2$.
We say that a path $\alpha:[0,1]\to D_2$ is
{\it in general position with respect to $\Gamma$}
if it avoids the vertices of $\Gamma$
and every intersection of $\alpha$ and $\Gamma$
is a transverse double point.
If $p$ is an intersection of $\alpha$ and an edge of $\Gamma$
labelled $i$ and 
if the edge is oriented from right to left 
(resp. from left to right),
then assign the intersection $p$ a letter $\sigma_i$
(resp. $\sigma_i^{-1}$).
Read the letters assigned to the intersections
of $\alpha$ and $\Gamma$ along $\alpha$ and
we have a word
$$\sigma_{i_1}^{\varepsilon_1}\sigma_{i_2}^{\varepsilon_2}\cdots \sigma_{i_s}^{\varepsilon_s}$$
in the braid generator.
We call this the {\it intersection braid word}
of $\alpha$ with respect to $\Gamma$, and denote it 
by $w_\Gamma(\alpha)$.

\begin{proposition}
{\rm (\cite[Proposition~10.4(1)]{QuandleInvariant})}
\label{DescribedByGamma}
Let $\Gamma$ be an $n$-chart in a $2$-disk $D_2$,
and $N(\Gamma)$ a regular neighborhood of $\Gamma$ 
in $D_2$.
Then there exists a simple surface braid $S$ in a $4$-dimensional disk $D_1\times D_2$ such that
\begin{enumerate}
\item[{\rm (a)}] the projection $pr_1:D_1\times D_2\to D_1$ satisfies 
the condition $pr_1(S\cap (D_1\times \{y\}))=X_n$ 
for any point $y$ in $Cl(D_2-N(\Gamma))$, where
$X_n$ denotes the $n$ fixed interior points in $D_1$,
\item[{\rm (b)}] the branched point set of $S$ corresponds 
to the set of the black vertices of $\Gamma$, and
\item[{\rm (c)}] for any simple path $\alpha:[0,1]\to D_2$
which is in general position with respect to $\Gamma$
and $\alpha(0),\alpha(1)$ are in $Cl(D_2 -N(\Gamma))$,
the set $S\cap (D_1\times \alpha([0,1]))$ is an $n$-braid
in $D_1\times \alpha([0,1])$ presented by 
the intersection braid word $w_\Gamma(\alpha)$.
\end{enumerate}
\end{proposition}

We say the surface $S$ in Proposition~\ref{DescribedByGamma}
is a {\it surface braid described by the chart $\Gamma$}.
We denote by $F(\Gamma)$
the closure of a surface braid described by the chart $\Gamma$.

\section{Knot quandles}
\label{s:KnotQuandle}

In this section, we review knot quandles and a relation between
quandle colorings and knot quandles.
We shall calculate knot quandles of surface braids described 
by the charts $T_k$ and $T_k^*$.

Let $D$ be a surface diagram of an oriented surface-link $F$.
Let ${\mathcal R}=\{R_1,\cdots,R_n\}$ be the set of all broken sheets of $D$.
Then we define the quandle $Q(D)$ generated by 
the set ${\mathcal R}$ and relations are as follows:
At each double point,
there exist three broken sheets $R_i,R_j,R_k\in {\mathcal R}$ uniquely such that 
$$R_i*R_j=R_k$$
where
$R_j$ is the over-sheet at the double point,
$R_i,R_k$ are under-sheets at the double point,
the normal of $R_j$ points from $R_i$ to $R_k$
(see Fig.~\ref{Fig09}).
Note that the quandle $Q(D)$ is isomorphic to the knot quandle $Q(F)$ of the surface-link $F$
(see \cite[Definition 10.7]{QuandleInvariant}).
Similarly we define the quandle $Q(S)$ of a surface braid $S$.

Let $Q$ be a quandle.
Let $D$ be a surface diagram of an oriented surface-link $F$,
and ${\mathcal R}=\{R_1,\cdots,R_n\}$ the set of all broken sheets of $D$.
For a quandle coloring ${\mathcal C}:{\mathcal R}\to Q$,
there is a homomorphism $\rho:Q(D)\to Q$
defined by $\rho(R_i)={\mathcal C}(R_i)$ for each broken sheet $R_i$.
Conversely,
for a homomorphism $\rho:Q(D)\to Q$,
there is a quandle coloring  ${\mathcal C}:{\mathcal R}\to Q$
defined by  ${\mathcal C}(R_i)=\rho(R_i)$ for each broken sheet $R_i$. Thus we have the following remark:

\begin{remark}
\label{RemarkColorHom}
Let $Q$ be a quandle, $F$ a surface-link or a surface braid,
and $Q(F)$ the knot quandle of $F$.
Then 
the set ${\rm Col}_Q(F)$ of the colorings
is naturally one-to-one corresponding
to the set ${\rm Hom}(Q(F),Q)$ of the homomorphisms
from $Q(F)$ to $Q$.
\end{remark}

Let $\Gamma$ be an $n$-chart in a 2-disk $D_2$,
and $y_1,y_2,\cdots,y_m$ all the black vertices of $\Gamma$.
Fix a point $y_0\in\partial D_2$.
Take a regular neighborhood $N(y_i)$ of $y_i$ in $D_2$
for each $i$.
Let $\alpha_i$ ($i=1,2,\cdots,m$) be a simple path
$\alpha_i:[0,1]\to Cl(D_2-(N(y_1)\cup\cdots\cup N(y_m)))$ 
such that
\begin{enumerate}
\item[(i)] $\alpha_i(0)\in \partial N(y_i)$ and $\alpha_i([0,1])\cap\partial D_2=\alpha_i(1)=y_0$ for each $i$,
\item[(ii)] $\alpha_i([0,1])\cap\alpha_j([0,1])=y_0$ for $i\not=j$,
\item[(iii)] the images of $\alpha_1,\alpha_2,\cdots,\alpha_m$ appear in this order around the point $y_0$.
\end{enumerate}
Then an $m$-tuple $(\alpha_1,\alpha_2,\cdots,\alpha_m)$ is called  a {\it Hurwitz arc system} for $\{y_1,y_2,\cdots,y_m\}$ with base point $y_0$.

\begin{proposition}
\label{LemmaQuandleInvariant}
{\rm (\cite[Lemma~10.9]{QuandleInvariant})}
Let $\Gamma$ be an $n$-chart in a $2$-disk $D_2$, and
$y_1,y_2,\cdots,y_m$ all the black vertices of $\Gamma$
contained in edges of label $k_1,k_2,\cdots,k_m$,
respectively.
Let $S$ be a surface braid described by $\Gamma$.
Let $(\alpha_1,\alpha_2,\cdots,\alpha_m)$ be a Hurwitz arc system for $\{y_1,y_2,\cdots,y_m\}$
such that each $\alpha_i$ is in general position with respect to $\Gamma$.
Then the quandle $Q(S)$ of the surface braid $S$ has a presentation whose generators are $x_1,x_2,\cdots,x_n$ and the relations are
$$Q(w_i)(x_{k_i})=Q(w_i)(x_{k_i +1}) \ \ (i=1,2,\cdots,m)$$
where $w_i=w_\Gamma(\alpha_i)$.
\end{proposition}

\begin{remark}
\label{RemarkKnotQuandle}
Let $Q$ be a quandle, and $F_Q\langle x_1,x_2,\cdots,x_n\rangle$ the free quandle generated by $x_1,x_2,\cdots,x_n$. 
If a homomorphism $c:F_Q\langle x_1,x_2,\cdots,x_n\rangle \to Q$
satisfies  $c(Q(w_i)(x_{k_i}))=c(Q(w_i)(x_{k_i +1})) \ \ (i=1,2,\cdots,m)$,
then there exists a homomorphism $\mathcal{C}:Q(S)\to Q$
with $\mathcal{C}\circ \iota=c$
where $\iota:F_Q\langle x_1,x_2,\cdots,x_n\rangle\to Q(S)$ is a natural homomorphism defined by $\iota(x_i)=x_i$ for $i=1,2,\cdots,n$.
\end{remark}

\begin{remark}
\label{DoesNotDepentBasePoint}
Let $S$ be a surface braid, and $\widehat{S}$ 
the closure of $S$.
Then the knot quandle $Q(S)$ is isomorphic to the
knot quandle $Q(\widehat{S})$.
Hence to calculate the knot quandle,
we do not depend a choice of a base point $y_0$ in $D_2-\Gamma$.
\end{remark}

\begin{lemma}
\label{LemmaConstantMap}
Let $Q$ be a quandle, $S$ a surface braid or a surface-link, and 
$F_Q\langle x_1,x_2,\cdots,x_n\rangle$ the free quandle.
Then we have the following:
\begin{enumerate}
\item[{\rm (a)}]
If a homomorphism $c:F_Q\langle x_1,x_2,\cdots,x_n\rangle\to Q$ 
satisfies
 $$c(Q(b)(x_1))=c(Q(b)(x_2))=\cdots=c(Q(b)(x_n))$$
for some $n$-braid $b$ in $B_n$,
then the map $c$ is a constant map.
\item[{\rm (a)}]
For a knot quandle $Q(S)$ generated 
by $x_1,x_2,\cdots,x_n$,
if a homomorphism
$\mathcal{C}:Q(S)\to Q$
satisfies
 $$\mathcal{C}(Q(b)(x_1))=\mathcal{C}(Q(b)(x_2))=\cdots=\mathcal{C}(Q(b)(x_n))$$
for some $n$-braid $b$ in $B_n$,
then the map $\mathcal{C}$ is a constant map.
\end{enumerate}
\end{lemma}

\begin{Proof}
{\bf Statement (a).}
Let $s=c(Q(b)(x_1))$.
By Lemma~\ref{QNii+1Equal},
for any braid $b'$ in $B_n$
we have 
$$s=c(Q(b'b)(x_1))=c(Q(b'b)(x_2))=\cdots=c(Q(b'b)(x_n)).$$
Thus if $b'=b^{-1}$,
then for each $i=1,2,\cdots,n$ we have
$$s=c(Q(b^{-1}b)(x_i))=c(Q(1)(x_i))=c(x_i).$$
Therefore
$c$ is a constant map.
Hence Statement (a) holds.

{\bf Statement (b).}
Let  $\iota:F_Q\langle x_1,x_2,\cdots,x_n\rangle\to Q(S)$
be the homomorphism defined by $\iota(x_i)=x_i$ for $i=1,2,\cdots,n$.
Let $c=\mathcal{C}\circ \iota$.
Then for $i=1,2,\cdots,n$, we have 
$c(Q(b)(x_i))=\mathcal{C}\circ \iota(Q(b)(x_i))=\mathcal{C}(Q(b)(x_i))$.
Hence
the map $c$ satisfies the condition (a) of this lemma.
By Lemma~\ref{LemmaConstantMap}(a),
the map $c$ is constant, so is $\mathcal{C}$.
Thus Statement (b) holds.
\end{Proof}


By Proposition~\ref{LemmaQuandleInvariant} 
and Remark~\ref{DoesNotDepentBasePoint},
we have the following lemma:

\begin{lemma}
\label{KnotQuandleTK}
Let $T_k$ be the $4$-chart as shown in Fig.~\ref{Fig02},
and $S_k$ a surface braid described by the chart $T_k$.
Let $(\alpha_1,\alpha_2,\cdots,\alpha_8)$ be a Hurwitz arc
system as shown in Fig.~\ref{Fig02}.
Then we have the intersection braid words as follows:
$$\begin{array}{rcl}
w_1=w_{T_k}(\alpha_1)& = & 1,\\
w_2=w_{T_k}(\alpha_2)& = & \sigma_2^{-(k+1)}\sigma_1\sigma_2^{2},\\
w_3=w_{T_k}(\alpha_3)&= &\sigma_2^{-2}\sigma_1\sigma_2^{2},\\
w_4=w_{T_k}(\alpha_4)& = & \sigma_2^{-(k+1)}\sigma_3\sigma_1\sigma_2^{2},\\
w_5(=w_4)=w_{T_k}(\alpha_5)& = & \sigma_2^{-(k+1)}\sigma_3\sigma_1\sigma_2^{2},\\
w_6=w_{T_k}(\alpha_6)& = & \sigma_2^{-2}\sigma_3\sigma_2^{2},\\
w_7=w_{T_k}(\alpha_7)& = &\sigma_2^{-(k+1)}\sigma_3\sigma_2^{2},\\
w_8(=w_1)=w_{T_k}(\alpha_8)& = & 1.
\end{array}$$
Moreover,
the quandle $Q(S_k)$ of the surface braid $S_k$
has a presentation whose generators are $x_1,x_2,x_3,x_4$
and the relations are

$$\begin{array}{rclcrcl}
Q(w_1)(x_1) & = & Q(w_1)(x_2), & & Q(w_5)(x_1) &=& Q(w_5)(x_2),\\
Q(w_2)(x_1)& = & Q(w_2)(x_2), & & Q(w_6)(x_1) &= & Q(w_6)(x_2),\\
Q(w_3)(x_3)& = & Q(w_3)(x_4), & & Q(w_7)(x_3)& = & Q(w_7)(x_4),\\
Q(w_4)(x_3) &= &  Q(w_4)(x_4), & & Q(w_8)(x_3)& = & Q(w_8)(x_4).\\\end{array}$$
\end{lemma}

\begin{lemma}
\label{KnotQuandleTK*}
Let $T_k,T_k^*$ be the $4$-charts as shown in Fig.~\ref{Fig02}, Fig.~\ref{Fig03}, respectively.
Let $S_k,S_k^*$ be surface braids described by the charts $T_k,T_k^*$, respectively.
Then the knot quandle $Q(S_k)$ is isomorphic to 
the knot quandle $Q(S_k^*)$.
\end{lemma}

\begin{Proof}
Note that we obtain from the chart $T_k$ 
to the chart $T_k^*$ by
changing the labels $1,3$ into
$3,1$, respectively.

Let $(\alpha_1',\alpha_2',\cdots,\alpha_8')$ be a Hurwitz arc
system as shown in Fig.~\ref{Fig03}.
Then we have the intersection braid words as follows:
$$\begin{array}{rcl}
w_1'=w_{T_k^*}(\alpha_1')& = & 1,\\
w_2'=w_{T_k^*}(\alpha_2')& = & \sigma_2^{-(k+1)}\sigma_3\sigma_2^{2},\\
w_3'=w_{T_k^*}(\alpha_3')&= &\sigma_2^{-2}\sigma_3\sigma_2^{2},\\
w_4'=w_{T_k^*}(\alpha_4')& = & \sigma_2^{-(k+1)}\sigma_1\sigma_3\sigma_2^{2},\\
w_5'(=w_4')=w_{T_k^*}(\alpha_5')& = & \sigma_2^{-(k+1)}\sigma_1\sigma_3\sigma_2^{2},\\
w_6'=w_{T_k^*}(\alpha_6')& = & \sigma_2^{-2}\sigma_1\sigma_2^{2},\\
w_7'=w_{T_k^*}(\alpha_7')& = &\sigma_2^{-(k+1)}\sigma_1\sigma_2^{2},\\
w_8'(=w_1')=w_{T_k^*}(\alpha_8')& = & 1.
\end{array}$$

Since $\sigma_1\sigma_3=\sigma_3\sigma_1$,
we have 
$w_4'=w_5'=w_4=w_5$,
where $w_i=w_{T_k}(\alpha_i)$ for $i=1,2,\cdots,8$.
Moreover
$w_1'=w_8'=w_1=w_8$,
$w_2'=w_7$, $w_3'=w_6$,
 $w_6'=w_3$ and $w_7'=w_2$.

Furthermore, by Proposition~\ref{LemmaQuandleInvariant} 
and Remark~\ref{DoesNotDepentBasePoint},
the quandle $Q(S_k^*)$ of the surface braid $S_k^*$
has a presentation whose generators are $x_1,x_2,x_3,x_4$
and the relations are

$$\begin{array}{rclcrcl}
Q(w_1')(x_3) & = & Q(w_1')(x_4),& & Q(w_5')(x_3) &=& Q(w_5')(x_4),\\
Q(w_2')(x_3)& = & Q(w_2')(x_4),& & Q(w_6')(x_3) &= & Q(w_6')(x_4),\\
Q(w_3')(x_1)& = & Q(w_3')(x_2),& & Q(w_7')(x_1)& = & Q(w_7')(x_2),\\
Q(w_4')(x_1) &= &  Q(w_4')(x_2), & & Q(w_8')(x_1)& = & Q(w_8')(x_2).\\
\end{array}$$
Therefore
$$\begin{array}{rclcrcl}
Q(w_8)(x_3) & = & Q(w_8)(x_4),& & Q(w_4)(x_3) &=& Q(w_4)(x_4),\\
Q(w_7)(x_3)& = & Q(w_7)(x_4),& & Q(w_3)(x_3) &= & Q(w_3)(x_4),\\
Q(w_6)(x_1)& = & Q(w_6)(x_2),& & Q(w_2)(x_1)& = & Q(w_2)(x_2),\\
Q(w_5)(x_1) &= &  Q(w_5)(x_2),& & Q(w_1)(x_1)& = & Q(w_1)(x_2).\\
\end{array}$$
Hence the knot quandle $Q(S_k^*)$ has the same relations of the knot quandle $Q(S_k)$.
Thus  $Q(S_k^*)$ is isomorphic to $Q(S_k)$.
\end{Proof}

Similarly, we have the following lemma:

\begin{lemma}
\label{KnotQuandleTK0}
Let $T_0$ be the $4$-chart as shown in Fig.~\ref{Fig04},
and $S_0$ a surface braid described by the chart $T_0$.
Let $(\alpha_1,\alpha_2,\cdots,\alpha_8)$ be a Hurwitz arc
system as shown in Fig.~\ref{Fig04}.
Then we have the intersection braid words as follows:
$$\begin{array}{rcl}
w_1=w_{T_0}(\alpha_1)& = & 1,\\
w_2=w_{T_0}(\alpha_2)& = & \sigma_1\sigma_2^{-1},\\
w_3(=w_2)=w_{T_0}(\alpha_3)&= & \sigma_1\sigma_2^{-1},\\
w_4=w_{T_0}(\alpha_4)& = & \sigma_2^{-1}\sigma_3\sigma_1\sigma_2^{-1},\\
w_5(=w_4)=w_{T_0}(\alpha_5)& = & \sigma_2^{-1}\sigma_3\sigma_1\sigma_2^{-1},\\
w_6=w_{T_0}(\alpha_6)& = & \sigma_3\sigma_2^{-1},\\
w_7(=w_6)=w_{T_0}(\alpha_7)& = &\sigma_3\sigma_2^{-1},\\
w_8(=w_1)=w_{T_0}(\alpha_8)& = & 1.
\end{array}$$
Moreover,
the quandle $Q(S_0)$ of the surface braid $S_0$
has a presentation whose generators are $x_1,x_2,x_3,x_4$
and the relations are

$$\begin{array}{rclcrcl}
Q(w_1)(x_1) & = & Q(w_1)(x_2), & & Q(w_5)(x_1) &= & Q(w_5)(x_2),\\
Q(w_2)(x_2)& = & Q(w_2)(x_3),& & Q(w_6)(x_2)& = & Q(w_6)(x_3).\\
Q(w_4)(x_3) &= &  Q(w_4)(x_4),& & Q(w_8)(x_3)& = & Q(w_8)(x_4).\\
\end{array}$$
\end{lemma}

\begin{remark}
\label{Remarkx1=x2x3=x4}
In the knot quandles $Q(S_0),Q(S_k),Q(S_k^*)$,
we have
$$x_1=x_2, \ {\rm and} \ x_3=x_4.$$
Because, $w_1=w_8=1$ implies that
$x_1=Q(w_1)(x_1)=Q(w_1)(x_2)=x_2$ and
$x_3=Q(w_8)(x_3)=Q(w_8)(x_4)=x_4$.
\end{remark}

\section{Colorings of a surface braid described by the 4-chart $T_0$}
\label{s:T0}

Let $S_0$ be a surface braid described by the 4-chart $T_0$
as shown in Fig.~\ref{Fig04}.
Then we have the quandle $Q(S_0)$ of the surface braid $S_0$.
In this section we shall count the number of homomorphims from $Q(S_0)$ to the quandle $Q_N$ in Example~\ref{PseudoTrivialQuandle}.

\begin{lemma}
\label{T0SSTT}
Let $N$ be an integer with $N\geqq3$.
Let $s,t$ be elements in $Q_N$ different from $N$.
Then there exists a homomorphism
$\mathcal{C}:Q(S_0)\to Q_N$ with
 $\mathcal{C}(x_1)=\mathcal{C}(x_2)=s$ and $\mathcal{C}(x_3)=\mathcal{C}(x_4)=t$.
\end{lemma}

\begin{Proof}
Let $F_Q\langle x_1,x_2,x_3,x_4\rangle$ be the free quandle generatied by $x_1,x_2,x_3,x_4$.
Let $c:F_Q\langle x_1,x_2,x_3,x_4\rangle \to Q_N$
be the homomorphism defined by
 $c(x_1)=c(x_2)=s$ and $c(x_3)=c(x_4)=t$.

We shall show $c(Q(w_2)(x_2))=c(Q(w_2)(x_3))$,
where $w_2=\sigma_1\sigma_2^{-1}$.
By Lemma~\ref{QNii+1NotN}\vspace{2mm}

$\begin{array}{cl}
(1) & (c(Q(\sigma_2^{-1})(x_1)),c(Q(\sigma_2^{-1})(x_2)),c(Q(\sigma_2^{-1})(x_3)),c(Q(\sigma_2^{-1})(x_4))) \\
&=(c(Q(1)(x_1)),c(Q(1)(x_3)),c(Q(1)(x_2)),c(Q(1)(x_4))) \\
&=(c(x_1),c(x_3),c(x_2),c(x_4)) \\
& = (s,t,s,t).\\
(2) & (c(Q(\sigma_1\sigma_2^{-1})(x_1)),c(Q(\sigma_1\sigma_2^{-1})(x_2)),c(Q(\sigma_1\sigma_2^{-1})(x_3)),c(Q(\sigma_1\sigma_2^{-1})(x_4)))\\
& =  (t,s,s,t).
\end{array}$\vspace{2mm}\\
Since $w_2=\sigma_1\sigma_2^{-1}$,
we have $c(Q(w_2)(x_2))=c(Q(w_2)(x_3))$.

We shall show $c(Q(w_6)(x_2))=c(Q(w_6)(x_3))$
where $w_6=\sigma_3\sigma_2^{-1}$.
By (1) and by Lemma~\ref{QNii+1NotN}\vspace{2mm}\\
$\begin{array}{l}
(c(Q(\sigma_3\sigma_2^{-1})(x_1)),c(Q(\sigma_3\sigma_2^{-1})(x_2)),c(Q(\sigma_3\sigma_2^{-1})(x_3)),c(Q(\sigma_3\sigma_2^{-1})(x_4))\\
 =  (s,t,t,s).
\end{array}$\vspace{2mm}\\
Since $w_6=\sigma_3\sigma_2^{-1}$,
we have $c(Q(w_6)(x_2))=c(Q(w_6)(x_3))$.

We shall show $c(Q(w_4)(x_3))=c(Q(w_4)(x_4))$
and $c(Q(w_5)(x_1))=c(Q(w_5)(x_2))$
where $w_4=w_5=\sigma_2^{-1}\sigma_3\sigma_1\sigma_2^{-1}$.
By (2) and by Lemma~\ref{QNii+1NotN}\vspace{2mm}

$\begin{array}{l}
(c(Q(\sigma_3\sigma_1\sigma_2^{-1})(x_1)),c(Q(\sigma_3\sigma_1\sigma_2^{-1})(x_2)),c(Q(\sigma_3\sigma_1\sigma_2^{-1})(x_3)),c(Q(\sigma_3\sigma_1\sigma_2^{-1})(x_4))\\
 =  (t,s,t,s).\\
(c(Q(w_4)(x_1)),c(Q(w_4)(x_2)),c(Q(w_4)(x_3)),c(Q(w_4)(x_4))\\
 =  (t,t,s,s).
\end{array}$\vspace{2mm}\\
Since $w_4=w_5$,
 we have $c(Q(w_4)(x_3))=c(Q(w_4)(x_4))$
and $c(Q(w_5)(x_1))=c(Q(w_5)(x_2))$.

Moreover, 
let $w_1=w_8=1$, then 
we have $c(Q(w_1)(x_1))=c(Q(w_1)(x_2))$ and 
$c(Q(w_8)(x_3))=c(Q(w_8)(x_4))$.
Therefore by Remark~\ref{RemarkKnotQuandle} and Lemma~\ref{KnotQuandleTK0},
we have the desired result.
\end{Proof}

\begin{lemma}
\label{T011NN}
Let $N$ be an integer with $N\geqq3$.
Then there does not exist a homomorphism
$\mathcal{C}:Q(S_0)\to Q_N$ with
 $\mathcal{C}(x_1)=\mathcal{C}(x_2)=1$ and $\mathcal{C}(x_3)=\mathcal{C}(x_4)=N$.
\end{lemma}

\begin{Proof}
Suppose that there exists a homomorphism
$\mathcal{C}:Q(S_0)\to Q_N$ with
 $\mathcal{C}(x_1)=\mathcal{C}(x_2)=1$ and $\mathcal{C}(x_3)=\mathcal{C}(x_4)=N$.
Let $\iota$ be the homomorphism from the
free quandle $F_Q\langle x_1,x_2,x_3,x_4\rangle$ to $Q(S_0)$
defined by $\iota(x_i)=x_i$ for $i=1,2,3,4$.
Then $\mathcal{C}\circ \iota$ is a homomorphism.
Let $c=\mathcal{C}\circ \iota$.
Then we have $c(x_1)=c(x_2)=1$ and $c(x_3)=c(x_4)=N.$

We shall show $c(Q(w_2)(x_2))\not=c(Q(w_2)(x_3))$,
where $w_2=\sigma_1\sigma_2^{-1}$.\vspace{2mm}

$\begin{array}{l}
 (c(Q(\sigma_2^{-1})(x_1)),c(Q(\sigma_2^{-1})(x_2)),c(Q(\sigma_2^{-1})(x_3)),c(Q(\sigma_2^{-1})(x_4))) \\
=(c(Q(1)(x_1)),c(Q(1)(x_3)),c(Q(1)(x_2)*Q(1)(x_3)),c(Q(1)(x_4))) \\
=(c(x_1),c(x_3),c(x_2* x_3),c(x_4)) \\
 = (1,N,1*N,N)=(1,N,2,N).\\
 (c(Q(\sigma_1\sigma_2^{-1})(x_1)),c(Q(\sigma_1\sigma_2^{-1})(x_2)),c(Q(\sigma_1\sigma_2^{-1})(x_3)),c(Q(\sigma_1\sigma_2^{-1})(x_4)))\\
=  (N\overline{*}1,1,2,N)=(N,1,2,N).
\end{array}$\vspace{2mm}\\
Since $w_2=\sigma_1\sigma_2^{-1}$,
we have $c(Q(w_2)(x_2))\not=c(Q(w_2)(x_3))$.

On the other hand, by Lemma~\ref{KnotQuandleTK0},
we have $Q(w_2)(x_2)=Q(w_2)(x_3)$ in $Q(S_0)$.
Thus the condition $c=\mathcal{C}\circ \iota$ implies
that $$c(Q(w_2)(x_2))=\mathcal{C}(Q(w_2)(x_2))=\mathcal{C}(Q(w_2)(x_3))=c(Q(w_2)(x_3)).$$
This is a contradiction.
Therefore there does not exist a homomorphism
$\mathcal{C}:Q(S_0)\to Q_N$ with
 $\mathcal{C}(x_1)=\mathcal{C}(x_2)=1$ and $\mathcal{C}(x_3)=\mathcal{C}(x_4)=N$.
\end{Proof}

\begin{lemma}
\label{T0NN11}
Let $N$ be an integer with $N\geqq3$.
Then there does not exist a homomorphism
$\mathcal{C}:Q(S_0)\to Q_N$ with
 $\mathcal{C}(x_1)=\mathcal{C}(x_2)=N$ and $\mathcal{C}(x_3)=\mathcal{C}(x_4)=1$.
\end{lemma}

\begin{Proof}
Suppose that there exists a homomorphism
$\mathcal{C}:Q(S_0)\to Q_N$ with
 $\mathcal{C}(x_1)=\mathcal{C}(x_2)=N$ and $\mathcal{C}(x_3)=\mathcal{C}(x_4)=1$.
Let $\iota$ be the homomorphism from the
free quandle $F_Q\langle x_1,x_2,x_3,x_4\rangle$ to $Q(S_0)$
defined by $\iota(x_i)=x_i$ for $i=1,2,3,4$.
Then $\mathcal{C}\circ \iota$ is a homomorphism.
Let $c=\mathcal{C}\circ \iota$.
Then we have
 $c(x_1)=c(x_2)=N$ and $c(x_3)=c(x_4)=1.$

We shall show $c(Q(w_6)(x_2))\not=c(Q(w_6)(x_3))$,
where $w_6=\sigma_3\sigma_2^{-1}$.\vspace{2mm}

$\begin{array}{l}
 (c(Q(\sigma_2^{-1})(x_1)),c(Q(\sigma_2^{-1})(x_2)),c(Q(\sigma_2^{-1})(x_3)),c(Q(\sigma_2^{-1})(x_4))) \\
 = (N,1,N*1,1)=(N,1,N,1).\\
(c(Q(\sigma_3\sigma_2^{-1})(x_1)),c(Q(\sigma_3\sigma_2^{-1})(x_2)),c(Q(\sigma_3\sigma_2^{-1})(x_3)),c(Q(\sigma_3\sigma_2^{-1})(x_4)))\\
=  (N,1,1\overline{*}N,N)=(N,1,N-1,N).
\end{array}$\vspace{2mm}\\
Since $w_6=\sigma_3\sigma_2^{-1}$,
we have $c(Q(w_6)(x_2))\not=c(Q(w_6)(x_3))$.

On the other hand, by Lemma~\ref{KnotQuandleTK0},
we have $Q(w_6)(x_2)=Q(w_6)(x_3)$ in $Q(S_0)$.
Thus by the similar way of the proof of Lemma~\ref{T011NN},
we have the same contradiction.
Therefore there does not exist a homomorphism
$\mathcal{C}:Q(S_0)\to Q_N$ with
 $\mathcal{C}(x_1)=\mathcal{C}(x_2)=N$ and $\mathcal{C}(x_3)=\mathcal{C}(x_4)=1$.
\end{Proof}

\begin{lemma}
\label{HomomorphismF}
Let $N$ be an integer with $N\geqq 3$.
Let $Q_N$ be the quandle in Example~\ref{PseudoTrivialQuandle}.
Let $f:Q_N\to Q_N$ be a map defined by

$f(x)=\left\{
\begin{array}{ll}
x+1, \ & ({\rm if} \ x<N-1),\\
1,\ & ({\rm if} \ x=N-1),\\
N,\ & ({\rm if} \ x=N).\\
\end{array}\right.$\\
Then $f$ is a homomorphism.
Moreover the inverse map $f^{-1}$ is also homomorphism
with 

$f^{-1}(x)=\left\{
\begin{array}{ll}
x-1, \ & ({\rm if} \ 1<x\leqq N-1),\\
N-1,\ & ({\rm if} \ x=1),\\
N,\ & ({\rm if} \ x=N).\\
\end{array}\right.$
\end{lemma}

\begin{Proof}
Let $x,y$ be elements in $Q_N$.
Then there are four cases:
(i) $y\not=N$,
(ii) $x<N-2, y=N$,
(iii) $x=N-1, y=N$,
(iv) $x=N, y=N$.

{\bf Case (i).}
Since $x*y=x$ and $f(y)\not=N$,
we have
$f(x)*f(y)=f(x)=f(x*y)$.

{\bf Case (ii).}
Since $f(x)=x+1<N-1$,
we have $f(x)*f(y)=(x+1)*N=x+2$.
Moreover, $x*y=x*N=x+1$ and
$f(x*y)=f(x+1)=x+2$.
Thus $f(x*y)=f(x)*f(y)$.

{\bf Case (iii).}
Since $x=N-1$, we have $f(x)=1$.
Thus the condition $y=N$ implies $f(x)*f(y)=1*N=2$.
Moreover, $f(x*y)=f((N-1)*N)=f(1)=2$.
Hence $f(x*y)=f(x)*f(y)$.

{\bf Case (vi).}
Since $x=N$ and $y=N$,
we have 
$f(x*y)=f(N*N)=f(N)=N$ and
$f(x)*f(y)=N*N=N$.
Thus $f(x*y)=f(x)*f(y)$.
\end{Proof}

From the above lemma,
we have the following lemma:

\begin{lemma}
\label{Homomorphism1ToSNtoN}
Let $N$ be an integer with $N\geqq 3$.
Let $Q_N$ be the quandle in Example~\ref{PseudoTrivialQuandle}.
Let $F$ be a surface-link or a surface braid, and 
$Q(F)$ the knot quandle of $F$.
Let $\mathcal{C}:Q(F)\to Q_N$ be a homomorphism
and $s$ an integer with $1\leqq s\leqq N-2$.
Then 
$f^{s}\circ \mathcal{C}$ and $f^{-s}\circ \mathcal{C}$ 
are homomorphisms from $Q(F)$ to $Q_N$,
where $f$ is the homomorphism in Lemma~\ref{HomomorphismF}.
Moreover, for an element $x$ in $Q(F)$,
\begin{enumerate}
\item[{\rm (a)}]
if $\mathcal{C}(x)=1$, then $f^{s}\circ \mathcal{C}(x)=s+1$, 
\item[{\rm (b)}]
if $\mathcal{C}(x)=s+1$, then $f^{-s}\circ \mathcal{C}(x)=1$,
\item[{\rm (c)}]
 if $\mathcal{C}(x)=N$, then $f^{s}\circ \mathcal{C}(x)=N$
and $f^{-s}\circ \mathcal{C}(x)=N$.
\end{enumerate}
\end{lemma}

\begin{corollary}
\label{CorT011NNNN11}
Let $N$ be an integer with $N\geqq3$.
Let $s$ be an element in the quandle $Q_N$ with $s\not=N$.
Then \begin{enumerate}
\item[{\rm (a)}]
there does not exist a homomorphism
$\mathcal{C}:Q(S_0)\to Q_N$
with $\mathcal{C}(x_1)=\mathcal{C}(x_2)=s$
and $\mathcal{C}(x_3)=\mathcal{C}(x_4)=N$, and
\item[{\rm (b)}]
there does not exist a homomorphism
$\mathcal{C}:Q(S_0)\to Q_N$
with $\mathcal{C}(x_1)=\mathcal{C}(x_2)=N$
and $\mathcal{C}(x_3)=\mathcal{C}(x_4)=s$.
\end{enumerate}
\end{corollary}

\begin{Proof}
By Lemma~\ref{T011NN} and Lemma~\ref{T0NN11},
we can assume $s\not=1$.
Now suppose that there exists
 a homomorphism $\mathcal{C}:Q(S_0)\to Q_N$
with 
 $\mathcal{C}(x_1)=\mathcal{C}(x_2)=s$ and
$\mathcal{C}(x_3)=\mathcal{C}(x_4)=N$.
By Lemma~\ref{Homomorphism1ToSNtoN}(b),(c),
the map $f^{-(s-1)}\circ \mathcal{C}:Q(S_0)\to Q_N$
is a homomorphism with
$f^{-(s-1)}\circ \mathcal{C}(x_1)=f^{-(s-1)}\circ \mathcal{C}(x_2)=1$ and 
$f^{-(s-1)}\circ \mathcal{C}(x_3)=f^{-(s-1)}\circ \mathcal{C}(x_4)=N$,
where $f$ is the map in Lemma~\ref{HomomorphismF}.
However by Lemma~\ref{T011NN}
there does not exist such a map $f^{-(s-1)}\circ \mathcal{C}$.
This is a contradiction. 
Therefore such a map $\mathcal{C}$ does not exist.
Hence Statement~(a) holds.

By the similar way as above,
we can show Statement~(b).
\end{Proof}

\begin{proposition}
\label{COLORT0}
Let $N$ be an integer with $N\geqq3$.
Let $Q_N$ be the quandle in Example~\ref{PseudoTrivialQuandle}.
Then 
$$|{\rm Col}_{Q_N}(S_0)|=|{\rm Col}_{Q_N}(F(T_0))|=(N-1)^2+1.$$
\end{proposition}

\begin{Proof}
By Remark~\ref{Remarkx1=x2x3=x4},
we have $x_1=x_2$ and $x_3=x_4$ in $Q(S_0)$.
Thus if $\mathcal{C}:Q(S_0)\to Q_N$ 
is a homomorphism,
then 
 $\mathcal{C}(x_1)=\mathcal{C}(x_2)$ and $\mathcal{C}(x_3)=\mathcal{C}(x_4).$

Let $s,t$ be any elements in $Q_N$.
We must determine when does there exist
a homomorphism $\mathcal{C}:Q(S_0)\to Q_N$
with $\mathcal{C}(x_1)=\mathcal{C}(x_2)=t$ and 
$\mathcal{C}(x_3)=\mathcal{C}(x_4)=s.$
Thus there are three cases:
(i) $s\not=N$ and $t\not=N$,
(ii) $s=N$ and $t=N$,
(iii) $s\not=t$ and ($s=N$ or $t=N$).

{\bf Case (i).}
By Lemma~\ref{T0SSTT},
there exists a homomorphism $\mathcal{C}:Q(S_0)\to Q_N$
with 
$\mathcal{C}(x_1)=\mathcal{C}(x_2)=s$ and
$\mathcal{C}(x_3)=\mathcal{C}(x_4)=t$.

{\bf Case (ii).}
Since any constant map $\mathcal{C}:Q(S_0)\to Q_N$
is a homomorphism,
there exists a homomorphism $\mathcal{C}:Q(S_0)\to Q_N$ with 
$\mathcal{C}(x_1)=\mathcal{C}(x_2)=\mathcal{C}(x_3)=\mathcal{C}(x_4)=N$.

By Cases (i),(ii) and Remark~\ref{RemarkColorHom}, we have
$$|{\rm Col}_{Q_N}(S_0)|\geqq (N-1)^2+1.$$
Therefore, it suffices to prove that
Case (iii) does not occur.

{\bf Case (iii).}
By Corollary~\ref{CorT011NNNN11},
there does not exist
a homomorphism $\mathcal{C}:Q(S_0)\to Q_N$
with $\mathcal{C}(x_1)=\mathcal{C}(x_2)=t$ and 
$\mathcal{C}(x_3)=\mathcal{C}(x_4)=s.$
\end{Proof}

\section{Colorings of surface braids described by the 4-charts $T_k$}
\label{s:T2k}

Let $S_k$ be a surface braid described by the 4-chart $T_k$
as shown in Fig.~\ref{Fig02}.
In this section we shall count the number of homomorphims from 
the quandle $Q(S_{2k})$ to the quandle $Q_N$ in Example~\ref{PseudoTrivialQuandle}.

\begin{lemma}
\label{MathCal2N-5}
Let $N,k$ be integers with $N\geqq3$ and $k\geqq 0$.
Let $c:F_Q\langle x_1,x_2,x_3,x_4\rangle\to Q_N$ be a homomorphism.
Suppose that $c(x_1)=c(x_2)=1$ and $c(x_3)=c(x_4)=N$.
Then we have the following:
\begin{enumerate}
\item[{\rm (a)}]
 $c(Q(w_3)(x_3))=c(Q(w_3)(x_4))$
where $w_3=\sigma_2^{-2}\sigma_1\sigma_2^{2}$.
\item[{\rm (b)}]
If $k\equiv 2N-5({\rm mod}\ 2(N-1))$, then 
 $c(Q(w_2)(x_1))=c(Q(w_2)(x_2))$
where $w_2=\sigma_2^{-(k+1)}\sigma_1\sigma_2^{2}$.
\item[{\rm (c)}]
If $k\not\equiv 2N-5({\rm mod}\ 2(N-1))$, then 
 $c(Q(w_2)(x_1))\not=c(Q(w_2)(x_2))$.
\end{enumerate} 
\end{lemma}

\begin{Proof}
Since $(c(x_1),c(x_2),c(x_3),c(x_4))=(1,1,N,N)$,
we have \vspace{2mm}

$\begin{array}{l}
(c(Q(\sigma_2)(x_1)),c(Q(\sigma_2)(x_2)),c(Q(\sigma_2)(x_3)),
c(Q(\sigma_2)(x_4)))\\
 = (1,N\overline{*}1,1,N)=(1,N,1,N).\vspace{2mm}\\
\end{array}$\\
Hence\vspace{2mm}\\
$\begin{array}{rl} 
(1) &  (c(Q(\sigma_2^2)(x_1)),c(Q(\sigma_2^2)(x_2)),c(Q(\sigma_2^2)(x_3)),c(Q(\sigma_2^2)(x_4)))\\
 &  =  (1,1\overline{*}N,N,N)=(1,N-1,N,N).\vspace{2mm}\\
(2) & (c(Q(\sigma_1\sigma_2^{2})(x_1)),c(Q(\sigma_1\sigma_2^{2})(x_2)),c(Q(\sigma_1\sigma_2^{2})(x_3)),c(Q(\sigma_1\sigma_2^{2})(x_4)))\\
 & =  ((N-1)\overline{*}1,1,N,N)=(N-1,1,N,N).
\end{array}$\vspace{2mm}\\
Thus\vspace{2mm}\\
$\begin{array}{l} 
(c(Q(\sigma_2^{-1}\sigma_1\sigma_2^2)(x_1)),c(Q(\sigma_2^{-1}\sigma_1\sigma_2^2)(x_2)),c(Q(\sigma_2^{-1}\sigma_1\sigma_2^2)(x_3)),c(Q(\sigma_2^{-1}\sigma_1\sigma_2^2)(x_4)))\\
=  (N-1,N,1*N,N)=(N-1,N,2,N).\vspace{2mm}\\
 (c(Q(\sigma_2^{-2}\sigma_1\sigma_2^2)(x_1)),c(Q(\sigma_2^{-2}\sigma_1\sigma_2^2)(x_2)),c(Q(\sigma_2^{-2}\sigma_1\sigma_2^2)(x_3)),
c(Q(\sigma_2^{-2}\sigma_1\sigma_2^2)(x_4)))\\
=  (N-1,2,N*2,N)=(N-1,2,N,N).
\end{array}$\vspace{2mm}\\
Since $w_3=\sigma_2^{-2}\sigma_1\sigma_2^2$,
we have  
$c(Q(w_3)(x_3))=c(Q(w_3)(x_4)).$
Hence Statement~(a) holds.

We shall show Statements (b) and (c).
By (2), we have $c(Q(\sigma_1\sigma_2^2)(x_1))=N-1$.
Thus by the definition of $Q(b)$, we have
$$(3)\ \ \ 
c(Q(\sigma_2^{-(k+1)}\sigma_1\sigma_2^2)(x_1))=c(Q(\sigma_1\sigma_2^2)(x_1))=N-1.$$ 

By (2),
we have $c(Q(\sigma_1\sigma_2^2)(x_2))=1$ and
$c(Q(\sigma_1\sigma_2^2)(x_3))=N$.
If $k\equiv 2N-5({\rm mod}\ 2(N-1))$,
then there exists an integer $m$ with 
$k=(2N-5)+2m(N-1)$.
Thus $k+1$ is even and
$k+1=2((N-2)+m(N-1))$.
Hence by Corollary~\ref{CorQNiPi+1N}(a),
we have
$c(Q(\sigma_2^{-(k+1)}\sigma_1\sigma_2^2)(x_2))=N-1$.
Thus by (3),
the condition $w_2=\sigma_2^{-(k+1)}\sigma_1\sigma_2^2$
implies that $c(Q(w_2)(x_1))=c(Q(w_2)(x_2))$.
Hence Statement~(b) holds.

Similarly
if $k\not\equiv 2N-5({\rm mod}\ 2(N-1))$,
then there exist integers $s,m$ with $0\leqq s\leqq 2N-3$,
$s\not=2N-5$ and $k=s+2m(N-1)$.
Thus $k+1$ is either odd, or $k+1=2(s'+m(N-1))$
where $s=2s'-1$ and $s'\not=N-2$.
Thus by Corollary~\ref{CorQNiPi+1N}(a) and (c)
we have
$$c(Q(\sigma_2^{-(k+1)}\sigma_1\sigma_2^2)(x_2))\not=N-1.$$
Hence by (3) we have $c(Q(w_2)(x_1))\not=c(Q(w_2)(x_2))$.
Thus Statement (c) holds.
\end{Proof}

\begin{corollary}
\label{CorMathCal2N-5}
Let $N,k$ be integers with $N\geqq3$ and $k\geqq 0$.
Let $s$ be an element in the quandle $Q_N$ with $s\not=N$.
If $k\not\equiv 2N-5({\rm mod}\ 2(N-1))$, then 
there does not exist a homomorphism
$\mathcal{C}:Q(S_k)\to Q_N$
with $\mathcal{C}(x_1)=\mathcal{C}(x_2)=s$
and $\mathcal{C}(x_3)=\mathcal{C}(x_4)=N$.
\end{corollary}

\begin{Proof}
Suppose that $k\not\equiv 2N-5({\rm mod}\ 2(N-1))$.
First we show for the case $s=1$.
Suppose that 
there exists a homomorphism
$\mathcal{C}:Q(S_k)\to Q_N$
with $\mathcal{C}(x_1)=\mathcal{C}(x_2)=1$
and $\mathcal{C}(x_3)=\mathcal{C}(x_4)=N$.
Let $\iota$ be the homomorphism from the
free quandle $F_Q\langle x_1,x_2,x_3,x_4\rangle$ to $Q(S_{k})$
defined by $\iota(x_i)=x_i$ for $i=1,2,3,4$.
Then $\mathcal{C}\circ \iota$ is a homomorphism.
Let $c=\mathcal{C}\circ \iota$.
Then $c(x_1)=c(x_2)=1$ and $c(x_3)=c(x_4)=N$.

By Lemma~\ref{MathCal2N-5}(c),
the condition $k\not\equiv 2N-5({\rm mod}\ 2(N-1))$
implies
 $c(Q(w_2)(x_1))\not=c(Q(w_2)(x_2))$.

On the other hand, we have
$Q(w_2)(x_1)=Q(w_2)(x_2)$ in $Q(S_{k})$
by Lemma~\ref{KnotQuandleTK}.
Thus by the similar way of the proof of Lemma~\ref{T011NN},
we have the same contradiction.
Hence there does not exist a homomorphism
$\mathcal{C}:Q(S_k)\to Q_N$
with $\mathcal{C}(x_1)=\mathcal{C}(x_2)=1$
and $\mathcal{C}(x_3)=\mathcal{C}(x_4)=N$.

By the similar way of the proof of Corollary~\ref{CorT011NNNN11}(a), we can show that 
there does not exist a homomorphism
$\mathcal{C}:Q(S_k)\to Q_N$
with $\mathcal{C}(x_1)=\mathcal{C}(x_2)=s$
and $\mathcal{C}(x_3)=\mathcal{C}(x_4)=N$.
\end{Proof}

\begin{lemma}
\label{MathCal2N-52}
Let $N,k$ be integers with $N\geqq3$ and $k\geqq 0$.
Let $c:F_Q\langle x_1,x_2,x_3,x_4\rangle\to Q_N$ be a homomorphism.
Suppose that $c(x_1)=c(x_2)=N$ and $c(x_3)=c(x_4)=1$.
Then we have the following:
\begin{enumerate}
\item[{\rm (a)}]
 $c(Q(w_6)(x_1))=c(Q(w_6)(x_2))$
where $w_6=\sigma_2^{-2}\sigma_3\sigma_2^2$.
\item[{\rm (b)}]
If $k\equiv 2N-5({\rm mod}\ 2(N-1))$,  then 
 $c(Q(w_7)(x_3))=c(Q(w_7)(x_4))$
where $w_7=\sigma_2^{-(k+1)}\sigma_3\sigma_2^2$.
\item[{\rm (c)}]
If $k\not\equiv 2N-5({\rm mod}\ 2(N-1))$,  then 
 $c(Q(w_7)(x_3))\not=c(Q(w_7)(x_4))$.
\end{enumerate} 
\end{lemma}

\begin{Proof}
Since $(c(x_1),c(x_2),c(x_3),c(x_4))=(N,N,1,1)$,
we have \vspace{2mm}

$\begin{array}{l}
(c(Q(\sigma_2)(x_1)),c(Q(\sigma_2 )(x_2)),c(Q(\sigma_2 )(x_3)),c(Q(\sigma_2)(x_4)))\\
 = (N,1\overline{*}N,N,1)=(N,N-1,N,1).\vspace{2mm}\\
\end{array}$\\
Hence\vspace{2mm}\\
$\begin{array}{rl} 
(1) &  (c(Q(\sigma_2^2)(x_1)),c(Q(\sigma_2^2)(x_2)),c(Q(\sigma_2^2)(x_3)),c(Q(\sigma_2^2)(x_4)))\\
 &  =  (N,N\overline{*}(N-1),N-1,1)=(N,N,N-1,1).\vspace{2mm}\\
(2) & (c(Q(\sigma_3\sigma_2^2)(x_1)),c(Q(\sigma_3\sigma_2^2)(x_2)),c(Q(\sigma_3\sigma_2^2)(x_3)),c(Q(\sigma_3\sigma_2^2)(x_4)))\\
 & =  (N,N,1\overline{*}(N-1),N-1)=(N,N,1,N-1).
\end{array}$\vspace{2mm}\\
Thus \vspace{2mm}

$\begin{array}{l} 
(c(Q(\sigma_2^{-1}\sigma_3\sigma_2^2)(x_1)),c(Q(\sigma_2^{-1}\sigma_3\sigma_2^2)(x_2)),c(Q(\sigma_2^{-1}\sigma_3\sigma_2^2)(x_3)))\\
 =  (N,1,N*1)=(N,1,N).\\
(c(Q(\sigma_2^{-2}\sigma_3\sigma_2^2)(x_1)),c(Q(\sigma_2^{-2}\sigma_3\sigma_2^2)(x_2)),c(Q(\sigma_2^{-2}\sigma_3\sigma_2^2)(x_3)))\\
 =  (N,N,1*N)=(N,N,2).
\end{array}$\vspace{2mm}\\
Since $w_6=\sigma_2^{-2}\sigma_3\sigma_2^2$,
we have $c(Q(w_6)(x_1))=c(Q(w_6)(x_2))$.
Thus Statement~(a) holds.

We shall show Statements (b) and (c).
By (2), we have $c(Q(\sigma_3\sigma_2^2)(x_4))=N-1$.
Thus by the definition of $Q(b)$, we have
$$(3)\ \ \ 
c(Q(\sigma_2^{-(k+1)}\sigma_3\sigma_2^2)(x_4))=c(Q(\sigma_3\sigma_2^2)(x_4))=N-1.$$ 

By (2),
we have $c(Q(\sigma_3\sigma_2^2)(x_2))=N$ and
$c(Q(\sigma_3\sigma_2^2)(x_3))=1$.
If $k\equiv 2N-5({\rm mod}\ 2(N-1))$,
then there exists an integer $m$ with 
$k=(2N-5)+2m(N-1)$.
Thus $k+1$ is even and
$k+1=2((N-2)+m(N-1))$.
Hence by Corollary~\ref{CorQNiNi+1P}(b),
we have
$c(Q(\sigma_2^{-(k+1)}\sigma_3\sigma_2^2)(x_3))=N-1$.
Thus by (3),
the condition $w_7=\sigma_2^{-(k+1)}\sigma_3\sigma_2^2$ implies 
that $c(Q(w_7)(x_3))=c(Q(w_7)(x_4))$.
Hence Statement~(b) holds.

Similarly
if $k\not\equiv 2N-5({\rm mod}\ 2(N-1))$,
then there exist integers $s,m$ with $0\leqq s\leqq 2N-3$,
$s\not=2N-5$ and $k=s+2m(N-1)$.
Thus $k+1$ is either odd, or $k+1=2(s'+m(N-1))$
where $s=2s'-1$ and $s'\not=N-2$.
Thus by Corollary~\ref{CorQNiNi+1P}(b) and (d)
we have
$$c(Q(\sigma_2^{-(k+1)}\sigma_3\sigma_2^2)(x_3))\not=N-1.$$
Hence by (3) we have $c(Q(w_7)(x_3))\not=c(Q(w_7)(x_4))$.
Thus Statement~(c) holds.
\end{Proof}

By the similar way of the proof of 
Corollary~\ref{CorMathCal2N-5},
we have the following corollary:

\begin{corollary}
\label{CorMathCal2N-52}
Let $N,k$ be integers with $N\geqq3$ and $k\geqq 0$.
Let $s$ be an element in the quandle $Q_N$ with $s\not=N$.
If $k\not\equiv 2N-5({\rm mod}\ 2(N-1))$, then 
there does not exist a homomorphism
$\mathcal{C}:Q(S_k)\to Q_N$
with $\mathcal{C}(x_1)=\mathcal{C}(x_2)=N$
and $\mathcal{C}(x_3)=\mathcal{C}(x_4)=s$.
\end{corollary}

\begin{lemma}
\label{MathCal3}
Let $N,k$ be integers with $N\geqq3$ and $k\geqq 0$.
Let $c:F_Q\langle x_1,x_2,x_3,x_4\rangle\to Q_N$ be a homomorphism.
Let $s,t$ be elements in $Q_N$ different from $N$.
Suppose that $c(x_1)=c(x_2)=s$ and $c(x_3)=c(x_4)=t$.
Then we have the following:
\begin{enumerate}
\item[{\rm (a)}]
If $k$ is odd, then 
 $c(Q(w_2)(x_1))=c(Q(w_2)(x_2))$
where $w_2=\sigma_2^{-(k+1)}\sigma_1\sigma_2^2$.
\item[{\rm (b)}]
If $k$ is even and $s\not=t$, then 
 $c(Q(w_2)(x_1))\not=c(Q(w_2)(x_2))$.
\item[{\rm (c)}]
 $c(Q(w_3)(x_3))=c(Q(w_3)(x_4))$
where $w_3=\sigma_2^{-2}\sigma_1\sigma_2^2$.
\end{enumerate} 
\end{lemma}

\begin{Proof}
Since $(c(x_1),c(x_2),c(x_3),c(x_4))=(s,s,t,t)$,
by Lemma~\ref{QNii+1NotN} we have \vspace{2mm}

$\begin{array}{rl}
(1) & (c(Q(\sigma_2^2)(x_1)),c(Q(\sigma_2^2)(x_2)),c(Q(\sigma_2^2)(x_3)),c(Q(\sigma_2^2)(x_4))\\
 & =  (s,s,t,t).\\
\end{array}$\vspace{2mm}\\
Hence\vspace{2mm}\\
$\begin{array}{rl} 
(2) & (c(Q(\sigma_1\sigma_2^2)(x_1)),c(Q(\sigma_1\sigma_2^2)(x_2)),c(Q(\sigma_1\sigma_2^2)(x_3)),c(Q(\sigma_1\sigma_2^2)(x_4))\\
 & =  (s,s,t,t).
\end{array}$\vspace{2mm}\\
Since $w_2=\sigma_2^{-(k+1)}\sigma_1\sigma_2^2$,
 by Lemma~\ref{QNi1i+12} we have that 
if $k$ is odd, then 
$$(3)\ \ \ (c(Q(w_2)(x_1)),c(Q(w_2)(x_2)),c(Q(w_2)(x_3)),c(Q(w_2)(x_4)) =  (s,s,t,t),$$ and
if $k$ is even, then 
$$(c(Q(w_2)(x_1)),c(Q(w_2)(x_2)),c(Q(w_2)(x_3)),c(Q(w_2)(x_4)) =  (s,t,s,t).$$
Thus Statements (a) and (b) hold.

Statement (c) follows from (3) for the case $k=1$.
\end{Proof}

By the similar way of the proof of 
Corollary~\ref{CorMathCal2N-5},
we have the following corollary:

\begin{corollary}
\label{CorMathCal2N-53}
Let $N,k$ be integers with $N\geqq3$ and $k\geqq 0$.
Let $s,t$ be elements in the quandle $Q_N$ different from $N$.
If $k$ is even and $s\not=t$, then 
there does not exist a homomorphism
$\mathcal{C}:Q(S_k)\to Q_N$
with $\mathcal{C}(x_1)=\mathcal{C}(x_2)=s$
and $\mathcal{C}(x_3)=\mathcal{C}(x_4)=t$.
\end{corollary}

\begin{proposition}
\label{PropS2k}
Let $N,k$ be integers with $N\geqq3$ and $k\geqq0$.
If $\mathcal{C}:Q(S_{2k})\to Q_N$ is a homomorphism,
then $\mathcal{C}$ is a constant map.
Moreover
we have $|{\rm Col}_{Q_N}(S_{2k})|=|{\rm Col}_{Q_N}(F(T_{2k}))|=N$.
\end{proposition}

\begin{Proof}
By Remark~\ref{Remarkx1=x2x3=x4}, 
we have $x_1=x_2$ and $x_3=x_4$ in $Q(S_{2k})$.
Since $\mathcal{C}:Q(S_{2k})\to Q_N$ is a map,

(1) $\mathcal{C}(x_1)=\mathcal{C}(x_2)$ and $\mathcal{C}(x_3)=\mathcal{C}(x_4).$

Let $s=\mathcal{C}(x_1)$ and $t=\mathcal{C}(x_3)$.
Next we shall show $s=t$.

Suppose that $s\not=t$.
There are two cases:
(i) $s\not=N$ and $t\not=N$,
(ii) $s=N$ or $t=N$.

{\bf Case (i).}
Since $2k$ is even,
by Corollary~\ref{CorMathCal2N-53}
there does not exist such a map $\mathcal{C}$.
Thus Case (i) does not occur.

{\bf Case (ii).}
Since $2k$ is even,
we have $2k\not\equiv 2N-5({\rm mod}\ 2(N-1))$.
Hence by Corollary~\ref{CorMathCal2N-5} and
Corollary~\ref{CorMathCal2N-52}
there does not exist such a map $\mathcal{C}$.
Thus Case (ii) does not occur.

Therefore all the cases do not occur.
Thus $s=t$.
Finally, by Lemma~\ref{LemmaConstantMap}(b)
the map $\mathcal{C}$
is a constant map
(i.e. any quandle coloring is a trivial coloring).
\end{Proof}

\section{Colorings of a surface braid $S_{2k-1}$}
\label{s:T2k-1}

In this section we shall count the number of homomorphisms from $Q(S_{2k-1})$ to the quandle $Q_N$ in Example~\ref{PseudoTrivialQuandle}.
We shall show the main theorem (Theorem~\ref{MainTheorem}).

\begin{lemma}
\label{k=2N-511NN}
Let $k$ be a positive integer,
and $Q_{k+2}$ the quandle in Example~\ref{PseudoTrivialQuandle}.
Then there exists a homomorphism
$\mathcal{C}:Q(S_{2k-1})\to Q_{k+2}$ with
 $\mathcal{C}(x_1)=\mathcal{C}(x_2)=1$ and $\mathcal{C}(x_3)=\mathcal{C}(x_4)=k+2$.
\end{lemma}

\begin{Proof}
Let $c:F_Q\langle x_1,x_2,x_3,x_4\rangle \to Q_{k+2}$
be the homomorphism defined by
 $c(x_1)=c(x_2)=1$ and $c(x_3)=c(x_4)=k+2$.
Let $N=k+2$.
Then $c:F_Q\langle x_1,x_2,x_3,x_4\rangle \to Q_N$,
and\vspace{2mm} 

(1) $2N-5=2(k+2)-5=2k-1$.

We shall show $c(Q(w_4)(x_3))=c(Q(w_4)(x_4))$
and $c(Q(w_5)(x_1))=c(Q(w_5)(x_2))$,
where $w_4=w_5=\sigma_2^{-((2k-1)+1)}\sigma_3\sigma_1\sigma_2^2$.
By (2) in the proof of Lemma~\ref{MathCal2N-5},\vspace{2mm}\\
$\begin{array}{l}
(c(Q(\sigma_1\sigma_2^2)(x_1)),c(Q(\sigma_1\sigma_2^2)(x_2)),c(Q(\sigma_1\sigma_2^2)(x_3)),c(Q(\sigma_1\sigma_2^2)(x_4)))\\ 
= (N-1,1,N,N).\end{array}$\\
Thus\vspace{2mm}\\
$\begin{array}{l} 
 (c(Q(\sigma_3\sigma_1\sigma_2^2)(x_1)),c(Q(\sigma_3\sigma_1\sigma_2^2)(x_2)),c(Q(\sigma_3\sigma_1\sigma_2^2)(x_3)),c(Q(\sigma_3\sigma_1\sigma_2^2)(x_4)))\\
=  (N-1,1,N\overline{*}N,N)=(N-1,1,N,N).\vspace{2mm}\\
\end{array}$\vspace{2mm}\\
Since $2k-1=2N-5$ (i.e. $2k=2(N-2)$) by (1),
  Corollary~\ref{CorQNiPi+1N}(a) and (b) implies that 
$$
c(Q(\sigma_2^{-2k}\sigma_3\sigma_1\sigma_2^2)(x_2))=N-1, \  c(Q(\sigma_2^{-2k}\sigma_3\sigma_1\sigma_2^2)(x_3))=N.$$
Moreover by the definition of $Q(b)$, we have
 
$\begin{array}{l}
c(Q(\sigma_2^{-2k}\sigma_3\sigma_1\sigma_2^2)(x_1))=c(Q(\sigma_3\sigma_1\sigma_2^2)(x_1))=N-1,\\
c(Q(\sigma_2^{-2k}\sigma_3\sigma_1\sigma_2^2)(x_4))=c(Q(\sigma_3\sigma_1\sigma_2^2)(x_4))=N.
\end{array}
$\\ 
Since $w_4=w_5=\sigma_2^{-2k}\sigma_3\sigma_1\sigma_2^2$,
we have $c(Q(w_4)(x_3))=c(Q(w_4)(x_4))$
and $c(Q(w_5)(x_1))=c(Q(w_5)(x_2))$.

We shall show $c(Q(w_6)(x_1))=c(Q(w_6)(x_2))$,
where $w_6=\sigma_2^{-2}\sigma_3\sigma_2^2$.
By (1) in the proof of Lemma~\ref{MathCal2N-5}, we have 
$$(c(Q(\sigma_2^2)(x_1)),c(Q(\sigma_2^2)(x_2)),c(Q(\sigma_2^2)(x_3)),c(Q(\sigma_2^2)(x_4)))=(1,N-1,N,N).$$
Hence we have \vspace{2mm}

$\begin{array}{rl}
(2) & (c(Q(\sigma_3\sigma_2^2)(x_2)),c(Q(\sigma_3\sigma_2^2)(x_3)),c(Q(\sigma_3\sigma_2^2)(x_4)))\\
 & =  (N-1,N\overline{*}N,N)=(N-1,N,N).\vspace{2mm}\\
& (c(Q(\sigma_2^{-1}\sigma_3\sigma_2^2)(x_2)),c(Q(\sigma_2^{-1}\sigma_3\sigma_2^2)(x_3)),c(Q(\sigma_2^{-1}\sigma_3\sigma_2^2)(x_4)))\\
 & =  (N,(N-1)*N,N)=(N,1,N).\vspace{2mm}\\
(3) & (c(Q(\sigma_2^{-2}\sigma_3\sigma_2^2)(x_2)),c(Q(\sigma_2^{-2}\sigma_3\sigma_2^2)(x_3)),c(Q(\sigma_2^{-2}\sigma_3\sigma_2^2)(x_4)))\\
& =  (1,N*1,N)=(1,N,N).
\end{array}$\vspace{2mm}\\
Thus $c(Q(w_6)(x_2))=1$.

On the other hand, by the definition of $Q(b)$, 
we have $c(Q(w_6)(x_1))=c(Q(\sigma_2^{-2}\sigma_3\sigma_2^2)(x_1))=c(Q(1)(x_1))=c(x_1)=1$.
Hence we have 
$c(Q(w_6)(x_1))=c(Q(w_6)(x_2))$.

We shall show 
$c(Q(w_7)(x_3))=c(Q(w_7)(x_4))$
where $w_7=\sigma_2^{-((2k-1)+1)}\sigma_3\sigma_2^2$.
Since $c(Q(\sigma_3\sigma_2^{2})(x_4))=N$ by (2),
we have $c(Q(w_7)(x_4))=c(Q(\sigma_3\sigma_2^2)(x_4))=N$.

Since  $2k=2N-4=2(N-3)+2$ by (1) and since
$c(Q(\sigma_2^{-2}\sigma_3\sigma_2^2)(x_2))=1$ and $c(Q(\sigma_2^{-2}\sigma_3\sigma_2^2)(x_3))=N$ by (3), 
we have 
$$c(Q(w_7)(x_3))=c(Q(\sigma_2^{-2k}\sigma_3\sigma_2^2)(x_3))=c(Q(\sigma_2^{-2(N-3)}\sigma_2^{-2}\sigma_3\sigma_2^2)(x_3))=N$$
by Corollary~\ref{CorQNiPi+1N}(b).
Hence we have $c(Q(w_7)(x_3))=c(Q(w_7)(x_4))$.

Moreover,
since $2k-1=2N-5$ by (1),
we have $c(Q(w_2)(x_1))=c(Q(w_2)(x_2))$ and 
$c(Q(w_3)(x_3))=c(Q(w_3)(x_4))$
by Lemma~\ref{MathCal2N-5}(a) and (b)
where $w_2=\sigma_2^{-((2k-1)+1)}\sigma_1\sigma_2^2$
and $w_3=\sigma_2^{-2}\sigma_1\sigma_2^2$.
Furthermore, let $w_1=w_8=1$,
then $c(Q(w_1)(x_1))=c(Q(w_1)(x_2))$ and 
$c(Q(w_8)(x_3))=c(Q(w_8)(x_4))$.
Therefore by Remark~\ref{RemarkKnotQuandle} 
and Lemma~\ref{KnotQuandleTK},
we have the desired result.
\end{Proof}


\begin{lemma}
\label{k=2N-5NN11}
Let $k$ be a positive integer,
and $Q_{k+2}$ the quandle in Example~\ref{PseudoTrivialQuandle}.
Then there exists a homomorphism
$\mathcal{C}:Q(S_{2k-1})\to Q_{k+2}$ with
 $\mathcal{C}(x_1)=\mathcal{C}(x_2)=k+2$ and $\mathcal{C}(x_3)=\mathcal{C}(x_4)=1$.
\end{lemma}

\begin{Proof}
Let $c:F_Q\langle x_1,x_2,x_3,x_4\rangle \to Q_{k+2}$
be the homomorphism defined by
 $c(x_1)=c(x_2)=k+2$ and $c(x_3)=c(x_4)=1$.
Let $N=k+2$.
Then $c:F_Q\langle x_1,x_2,x_3,x_4\rangle \to Q_N$,
and\vspace{2mm} 

(1) $2N-5=2(k+2)-5=2k-1$.

We shall show $c(Q(w_4)(x_3))=c(Q(w_4)(x_4))$
and $c(Q(w_5)(x_1))=c(Q(w_5)(x_2))$,
where $w_4=w_5=\sigma_2^{-((2k-1)+1)}\sigma_3\sigma_1\sigma_2^2$.

By (1) in the proof of Lemma~\ref{MathCal2N-52}, \vspace{2mm}\\
$\begin{array}{l}
(c(Q(\sigma_2^2)(x_1)),c(Q(\sigma_2^2)(x_2)),c(Q(\sigma_2^2)(x_3)),c(Q(\sigma_2^2)(x_4))) 
= (N,N,N-1,1).\end{array}$\vspace{2mm}\\
Hence\vspace{2mm}\\
$\begin{array}{rl} 
(2) & (c(Q(\sigma_1\sigma_2^2)(x_1)),c(Q(\sigma_1\sigma_2^2)(x_2)),c(Q(\sigma_1\sigma_2^2)(x_3)),c(Q(\sigma_1\sigma_2^2)(x_4)))\\
&  =  (N\overline{*}N,N,N-1,1)=(N,N,N-1,1).\\
 & (c(Q(\sigma_3\sigma_1\sigma_2^2)(x_1)),c(Q(\sigma_3\sigma_1\sigma_2^2)(x_2)),c(Q(\sigma_3\sigma_1\sigma_2^2)(x_3)),c(Q(\sigma_3\sigma_1\sigma_2^2)(x_4)))\\
 &=  (N,N,1\overline{*}(N-1),N-1)=(N,N,1,N-1).\vspace{2mm}\\
\end{array}$\vspace{2mm}\\
Since $2k-1=2N-5$ (i.e. $2k=2(N-2)$) by (1), 
 Corollary~\ref{CorQNiNi+1P}(a) and (b) implies that
we have
$$
c(Q(\sigma_2^{-2k}\sigma_3\sigma_1\sigma_2^2)(x_2))=N,\ c(Q(\sigma_2^{-2k}\sigma_3\sigma_1\sigma_2^2)(x_3))=N-1.$$
Moreover by the definition of $Q(b)$, we have
 
$\begin{array}{l}
c(Q(\sigma_2^{-2k}\sigma_3\sigma_1\sigma_2^2)(x_1))=c(Q(\sigma_3\sigma_1\sigma_2^2)(x_1))=N,\\
c(Q(\sigma_2^{-2k}\sigma_3\sigma_1\sigma_2^2)(x_4))=c(Q(\sigma_3\sigma_1\sigma_2^2)(x_4))=N-1.
\end{array}
$\\ 
Since $w_4=w_5=\sigma_2^{-2k}\sigma_3\sigma_1\sigma_2^2$,
we have $c(Q(w_4)(x_3))=c(Q(w_4)(x_4))$
and $c(Q(w_5)(x_1))=c(Q(w_5)(x_2))$.

We shall show $c(Q(w_3)(x_3))=c(Q(w_3)(x_4))$,
where $w_3=\sigma_2^{-2}\sigma_1\sigma_2^2$.
By (2), we have \vspace{2mm} 
 
$\begin{array}{rl}
 & (c(Q(\sigma_2^{-1}\sigma_1\sigma_2^2)(x_2)),c(Q(\sigma_2^{-1}\sigma_1\sigma_2^2)(x_3)),c(Q(\sigma_2^{-1}\sigma_1\sigma_2^2)(x_4)))\\
  & =  (N-1,N*(N-1),1)=(N-1,N,1).\vspace{2mm}\\
(3) & (c(Q(\sigma_2^{-2}\sigma_1\sigma_2^2)(x_2)),c(Q(\sigma_2^{-2}\sigma_1\sigma_2^2)(x_3)),c(Q(\sigma_2^{-2}\sigma_1\sigma_2^2)(x_4)))\\
& =  (N,(N-1)*N,1)=(N,1,1).
\end{array}$\vspace{2mm}\\
Hence we have 
$c(Q(w_3)(x_3))=c(Q(w_3)(x_4))$.

We shall show 
$c(Q(w_2)(x_1))=c(Q(w_2)(x_2))$
where $w_2=\sigma_2^{-((2k-1)+1)}\sigma_1\sigma_2^2$.
Since $c(Q(\sigma_1\sigma_2^2)(x_1))=N$ by (2),
we have $c(Q(w_2)(x_1))=c(Q(\sigma_1\sigma_2^2)(x_1))=N$.

Since  $2k=2N-4=2(N-3)+2$ by (1) and since $c(Q(\sigma_2^{-2}\sigma_1\sigma_2^2)(x_2))=N$ and $c(Q(\sigma_2^{-2}\sigma_1\sigma_2^2)(x_3))=1$ by (3), 
we have 
$$c(Q(w_2)(x_2))=c(Q(\sigma_2^{-2k}\sigma_1\sigma_2^2)(x_2))=c(Q(\sigma_2^{-2(N-3)}\sigma_2^{-2}\sigma_1\sigma_2^2)(x_2))=N$$
by Corollary~\ref{CorQNiNi+1P}(a).
Hence we have $c(Q(w_2)(x_1))=c(Q(w_2)(x_2))$.

Moreover, since $2k-1=2N-5$ by (1),
we have $c(Q(w_6)(x_1))=c(Q(w_6)(x_2))$ and 
$c(Q(w_7)(x_3))=c(Q(w_7)(x_4))$ by Lemma~\ref{MathCal2N-52}(a) and (b)
where $w_6=\sigma_2^{-2}\sigma_3\sigma_2^2$
and
$w_7=\sigma_2^{-((2k-1)+1)}\sigma_3\sigma_2^2$.
Furthermore, let $w_1=w_8=1$,
then $c(Q(w_1)(x_1))=c(Q(w_1)(x_2))$ and 
$c(Q(w_8)(x_3))=c(Q(w_8)(x_4))$.
Therefore by Remark~\ref{RemarkKnotQuandle} and Lemma~\ref{KnotQuandleTK},
we have the desired result.
\end{Proof}

By Lemma~\ref{Homomorphism1ToSNtoN}(a),(c),
Lemma~\ref{k=2N-511NN} and
Lemma~\ref{k=2N-5NN11},
we have the following corollary:

\begin{corollary}
\label{k=2N-5ssNNNNss}
Let $k$ be a positive integer,
and $Q_{k+2}$ the quandle in Example~\ref{PseudoTrivialQuandle}.
Let $s$ be an element in $Q_{k+2}$ with $s\not=k+2$.
Then 
\begin{enumerate}
\item[{\rm (a)}]
there exists a homomorphism
$\mathcal{C}:Q(S_{2k-1})\to Q_{k+2}$ with
 $\mathcal{C}(x_1)=\mathcal{C}(x_2)=s$ and $\mathcal{C}(x_3)=\mathcal{C}(x_4)=k+2$, and 
\item[{\rm (b)}]
 there exists a homomorphism
$\mathcal{C}:Q(S_{2k-1})\to Q_{k+2}$ with
 $\mathcal{C}(x_1)=\mathcal{C}(x_2)=k+2$ and $\mathcal{C}(x_3)=\mathcal{C}(x_4)=s$.
\end{enumerate}
\end{corollary}

\begin{lemma}
\label{KoddSSTT}
Let $N,k$ be integers with $N\geqq3$ and $k\geqq0$.
Let $s,t$ be elements in $Q_N$ different from $N$.
Then there exists a homomorphism
$\mathcal{C}:Q(S_{2k-1})\to Q_N$ with
 $\mathcal{C}(x_1)=\mathcal{C}(x_2)=s$ and $\mathcal{C}(x_3)=\mathcal{C}(x_4)=t$.
\end{lemma}

\begin{Proof}
Let $c:F_Q\langle x_1,x_2,x_3,x_4\rangle \to Q_N$
be the homomorphism defined by
 $c(x_1)=c(x_2)=s$ and $c(x_3)=c(x_4)=t$.

We shall show $c(Q(w_4)(x_3))=c(Q(w_4)(x_4))$
and $c(Q(w_5)(x_1))=c(Q(w_5)(x_2))$,
where $w_4=w_5=\sigma_2^{-((2k-1)+1)}\sigma_3\sigma_1\sigma_2^2$.
By (2) in the proof of Lemma~\ref{MathCal3},\vspace{2mm}

$\begin{array}{l}
(c(Q(\sigma_1\sigma_2^2)(x_1)),c(Q(\sigma_1\sigma_2^2)(x_2)),c(Q(\sigma_1\sigma_2^2)(x_3)),c(Q(\sigma_1\sigma_2^2)(x_4))) \\
= (s,s,t,t).\end{array}$\\
Thus by Lemma~\ref{QNii+1NotN}\vspace{2mm}\\
$\begin{array}{l} 
(c(Q(\sigma_3\sigma_1\sigma_2^2)(x_1)),c(Q(\sigma_3\sigma_1\sigma_2^2)(x_2)),c(Q(\sigma_3\sigma_1\sigma_2^2)(x_3)),c(Q(\sigma_3\sigma_1\sigma_2^2)(x_4)))\\
 =  (s,s,t,t).
\end{array}$\vspace{2mm}\\
Since $2k$ is even,
by Lemma~\ref{QNi1i+12}(a)
we have
$$
(c(Q(\sigma_2^{-2k}\sigma_3\sigma_1\sigma_2^2)(x_2)),c(Q(\sigma_2^{-2k}\sigma_3\sigma_1\sigma_2^2)(x_3)))
= (s,t).$$
Moreover by the definition of $Q(b)$, we have
 
$\begin{array}{l}
c(Q(\sigma_2^{-2k}\sigma_3\sigma_1\sigma_2^2)(x_1))=c(Q(\sigma_3\sigma_1\sigma_2^2)(x_1))=s,\\
c(Q(\sigma_2^{-2k}\sigma_3\sigma_1\sigma_2^2)(x_4))=c(Q(\sigma_3\sigma_1\sigma_2^2)(x_4))=t.
\end{array}
$\\ 
Since $w_4=w_5=\sigma_2^{-2k}\sigma_3\sigma_1\sigma_2^2$,
we have $c(Q(w_4)(x_3))=c(Q(w_4)(x_4))$
and $c(Q(w_5)(x_1))=c(Q(w_5)(x_2))$.

We shall show $c(Q(w_6)(x_1))=c(Q(w_6)(x_2))$
and $c(Q(w_7)(x_3))=c(Q(w_7)(x_4))$
where $w_6=\sigma_2^{-2}\sigma_3\sigma_2^2$ and
 $w_7=\sigma_2^{-((2k-1)+1)}\sigma_3\sigma_2^2$.
By Lemma~\ref{QNi1i+12}(a), we have 
$$(c(Q(\sigma_2^2)(x_1)),c(Q(\sigma_2^2)(x_2)),c(Q(\sigma_2^2)(x_3)),c(Q(\sigma_2^2)(x_4))=(s,s,t,t).$$
Hence by Lemma~\ref{QNii+1NotN} we have \vspace{2mm}

$\begin{array}{l}
(c(Q(\sigma_3\sigma_2^2)(x_1)),c(Q(\sigma_3\sigma_2^2)(x_2)),c(Q(\sigma_3\sigma_2^2)(x_3)),c(Q(\sigma_3\sigma_2^2)(x_4))\\
 =  (s,s,t,t).
\end{array}$\vspace{2mm}\\
Since $w_6=\sigma_2^{-2}\sigma_3\sigma_2^2$,
 by Lemma~\ref{QNi1i+12}(a)
we have $c(Q(w_6)(x_1))=c(Q(w_6)(x_2))=s$.
Since $w_7=\sigma_2^{-2k}\sigma_3\sigma_2^2$,
 by Lemma~\ref{QNi1i+12}(a)
we have $c(Q(w_7)(x_3))=c(Q(w_7)(x_4))=t$.

Moreover,
since $2k-1$ is odd,
by Lemma~\ref{MathCal3}(a) and (c)
we have $c(Q(w_2)(x_1))=c(Q(w_2)(x_2))$ and 
$c(Q(w_3)(x_3))=c(Q(w_3)(x_4))$ 
where  $w_2=\sigma_2^{-((2k-1)+1)}\sigma_1\sigma_2^2$
and $w_3=\sigma_2^{-2}\sigma_1\sigma_2^2$.
Furthermore, let $w_1=w_8=1$,
then $c(Q(w_1)(x_1))=c(Q(w_1)(x_2))$ and 
$c(Q(w_8)(x_3))=c(Q(w_8)(x_4))$.
Therefore by Remark~\ref{RemarkKnotQuandle} and Lemma~\ref{KnotQuandleTK},
we have the desired result.
\end{Proof}

\begin{proposition}
\label{PropS2k-1}
Let $k$ be a positive integer.
Then $$|{\rm Col}_{Q_{k+2}}(S_{2k-1})|=
|{\rm Col}_{Q_{k+2}}(F(T_{2k-1}))|=(k+2)^2.$$
\end{proposition}

\begin{Proof}
By Remark~\ref{Remarkx1=x2x3=x4},
we have $x_1=x_2$ and $x_3=x_4$ in $Q(S_{2k-1})$.
Thus if $\mathcal{C}:Q(S_{2k-1})\to Q_{k+2}$ is a homomorphism,
then 

(1) $\mathcal{C}(x_1)=\mathcal{C}(x_2)$ and $\mathcal{C}(x_3)=\mathcal{C}(x_4).$

Let $s,t$ be any elements in $Q_{k+2}$.
By (1) and Remark~\ref{RemarkColorHom},
it suffices to prove that
there exists a homomorphism
$\mathcal{C}:Q(S_{2k-1})\to Q_{k+2}$
with 
$\mathcal{C}(x_1)=\mathcal{C}(x_2)=s$ and
$\mathcal{C}(x_3)=\mathcal{C}(x_4)=t$.

Since any constant map $\mathcal{C}:Q(S_{2k-1})\to Q_{k+2}$
is a homomorphism,
we can assume that $s\not=t$.
There are two cases:
(i) $s\not=k+2$ and $t\not=k+2$,
(ii) $s=k+2$ or $t=k+2$.

{\bf Case (i).}
By Lemma~\ref{KoddSSTT},
there exists such a map $\mathcal{C}$.

{\bf Case (ii).}
By Corollary~\ref{k=2N-5ssNNNNss},
there exists such a map $\mathcal{C}$.
\end{Proof}

\begin{proposition}
\label{PropS2l-1}
Let $k,\ell$ be positive integers with $\ell<k$.
Then $$|{\rm Col}_{Q_{k+2}}(S_{2\ell-1})|=
|{\rm Col}_{Q_{k+2}}(F(T_{2\ell-1}))|=(k+1)^2+1.$$
\end{proposition}

\begin{Proof}
By Remark~\ref{Remarkx1=x2x3=x4},
we have $x_1=x_2$ and $x_3=x_4$ in $Q(S_{2\ell-1})$.
Thus if $\mathcal{C}:Q(S_{2\ell-1})\to Q_{k+2}$ 
is a homomorphism,
then 
 $\mathcal{C}(x_1)=\mathcal{C}(x_2)$ and $\mathcal{C}(x_3)=\mathcal{C}(x_4).$

Let $s,t$ be any elements in $Q_{k+2}$.
We must determine when does there exist
a homomorphism $\mathcal{C}:Q(S_{2\ell-1})\to Q_{k+2}$
with $\mathcal{C}(x_1)=\mathcal{C}(x_2)=t$ and 
$\mathcal{C}(x_3)=\mathcal{C}(x_4)=s.$
Thus there are three cases:
(i) $s\not=k+2$ and $t\not=k+2$,
(ii) $s=k+2$ and $t=k+2$,
(iii) $s\not=t$ and ($s=k+2$ or $t=k+2$).

{\bf Case (i).}
By Lemma~\ref{KoddSSTT},
there exists a homomorphism $\mathcal{C}:Q(S_{2\ell-1})\to Q_{k+2}$
with 
$\mathcal{C}(x_1)=\mathcal{C}(x_2)=s$ and
$\mathcal{C}(x_3)=\mathcal{C}(x_4)=t$.

{\bf Case (ii).}
Since any constant map $\mathcal{C}:Q(S_{2\ell-1})\to Q_{k+2}$
is a homomorphism,
there exists a homomorphism $\mathcal{C}:Q(S_{2\ell-1})\to Q_{k+2}$ with 
$\mathcal{C}(x_1)=\mathcal{C}(x_2)=\mathcal{C}(x_3)=\mathcal{C}(x_4)=k+2$.

By Cases (i),(ii) and Remark~\ref{RemarkColorHom}, 
we have
$$|{\rm Col}_{Q_{k+2}}(S_{2\ell-1})|\geqq (k+1)^2+1.$$
Therefore, it suffices to prove that
Case (iii) does not occur.

{\bf Case (iii).}
Let $N=k+2$. 
Since $\ell<k$, we have
$2N-5=2(k+2)-5=2k-1>2\ell-1$.
Thus $2\ell-1\not\equiv 2N-5 ({\rm mod}\ 2(N-1))$.
Hence by Corollary~\ref{CorMathCal2N-5} and Corollary~\ref{CorMathCal2N-52}
there does not exist a homomorphism $\mathcal{C}:Q(S_{2\ell-1})\to Q_{k+2}$ with 
$\mathcal{C}(x_1)=\mathcal{C}(x_2)=s$ 
and $\mathcal{C}(x_3)=\mathcal{C}(x_4)=t$.
\end{Proof}

Finally we shall show Theorem~\ref{MainTheorem}.

{\it Proof of Theorem~\ref{MainTheorem}.}

Statement~(a) follow from Proposition~\ref{COLORT0}.

By Lemma~\ref{KnotQuandleTK*},
we have
$|{\rm Col}_{Q_N}(F(T_{k}))|=|{\rm Col}_{Q_N}(F(T_{k}^*))|$
for any positive integers $k,N$ with $N\geqq 3$.
Hence we only calculate the number
$|{\rm Col}_{Q_N}(F(T_{k}))|$.

Statement~(b) follows from Proposition~\ref{PropS2k}.
Statement~(c) follows from Proposition~\ref{PropS2k-1}.
Statement~(d) follows from Proposition~\ref{PropS2l-1}.
Therefore we complete the proof of Theorem~\ref{MainTheorem}.
{\hfill {$\square$}\vspace{1.5em}}





\vspace{5mm}

\begin{minipage}{65mm}
{Teruo NAGASE
\\
{\small Tokai University \\
4-1-1 Kitakaname, Hiratuka \\
Kanagawa, 259-1292 Japan\\
\\
nagase@keyaki.cc.u-tokai.ac.jp
}}
\end{minipage}
\begin{minipage}{65mm}
{Akiko SHIMA 
\\
{\small Department of Mathematics, 
\\
Tokai University
\\
4-1-1 Kitakaname, Hiratuka \\
Kanagawa, 259-1292 Japan\\
shima@keyaki.cc.u-tokai.ac.jp
}}
\end{minipage}


\begin{thebibliography}{99}

\bibitem{QuandleInvariant}J.\ S.\ Carter, D.\ Jelsovsky, S.\ Kamada, L.\ Langford, and M.\ Saito,
{\it Quandle cohomology and state-sum invariants of knotted curves and surfaces}, Trans. Amer. Math. Soc. {\bf 355} (2003), 
3947--3989, MR1990571 (2005b:57048).

\bibitem{Surfaces4Space}J.\ S.\ Carter, S.\ Kamada and M.\ Saito, {"Surfaces in $4$-Space"}, Encyclopaedia of Mathematical Sciences  {\bf 142} Low-Dimensional Topology, III. Springer-Verlag, Berlin, (2004), MR2060067 (2005e:57065).


\bibitem{KnottedSurfaces}J.\ S.\ Carter and M.\ Saito, {"Knotted surfaces and their diagrams"}, Mathematical Surveys and Monographs, vol. 55, American Mathematical Society, Providence, RI, 1998, MR1487374 (98m:57027).



\bibitem{BraidThree} S.\ Kamada, {\it Surfaces in $R^4$ of braid index three are ribbon}, J. Knot Theory Ramif. {\bf 1}(2) (1992) 137--160, MR1164113 (93h:57039).

\bibitem{CharacterKamada} S.\ Kamada, {\it A characterization of grups of closed orientable surfaces in 4-space}, Topology {\bf 33} (1994) 113--122, MR1259518 (95a:57002).

\bibitem{BraidBook} S.\ Kamada, {"Braid and Knot Theory in Dimension Four"}, 
Mathematical Surveys and Monographs,
vol. 95, American Mathematical Society, 2002, MR1900979 (2003d:57050).

\bibitem{NagaseHirota} T.\ Nagase and A.\ Hirota, {\it The closure of a surface braid represented by a 4-chart with at most one crossing is a ribbon surface}, Osaka J. Math. {\bf 43} (2006) 413--430, MR2262343 (2007g:57040).


\bibitem{MinimalChartOneCrossing} T.\ Nagase and A.\ Shima, 
{\it Any chart with at most one crossing is a ribbon chart}, 
Topology Appl. {\bf 157} (2010), 1703--1720, MR2639836 (2011f:57048).

\bibitem{NSTwoCrossing} T.\ Nagase and A.\ Shima, {\it On charts with two crossings I: there exist no NS-tangles in a minimal chart}, J. Math. Sci. Univ. Tokyo {\bf 17} (2010) 217--241, MR2759760 (2012a:57032).


\bibitem{NSTwoCrossingII} T.\ Nagase and A.\ Shima, {\it On charts with two crossings II}, Osaka J. Math. {\bf 49} (2012) 909--929, MR3007949.




\bibitem{StI} T.\ Nagase and A.\ Shima, 
{\it The structure of a minimal $n$-chart with two crossings I: Complementary domains of $\Gamma_1\cup \Gamma_{n-1}$}, 
J. Knot Theory Ramif. {\bf 27}(14) (2018) 1850078, 37 pages, arXiv:1704.01232v3, MR3896318.

\bibitem{StII} T.\ Nagase and A.\ Shima, 
{\it The structure of a minimal $n$-chart with two crossings II: Neighbourhoods of $\Gamma_1\cup\Gamma_{n-1}$},
Revista de la Real Academia de Ciencias Exactas, Fiskcas y Natrales. Serie A. Mathem\'aticas {\bf 113} (2019) 1693--1738, arXiv:1709.08827v2, MR3956213.


\bibitem{NST} T.\ Nagase, A.\ Shima and H.\ Tsuji, 
{\it The closures of surface braids obtained from minimal $n$-charts with four white vertices}, 
J. Knot Theory Ramif. {\bf 22}(2) (2013) 1350007, 27 pages, MR3037298. 


\bibitem{Tanaka} K.\ Tanaka, 
{\it A Note on CI-moves}, 
Intelligence of low dimensional topology 2006,
Ser. Knots Everything, vol. 40, World Sci. Publ., Hackensack, NJ, 2007, pp.307--314, MR2371740 (2009a:57017).
\end{thebibliography}
\end{document}